\documentclass[12 pt]{amsart}
\usepackage[english]{babel}
\usepackage[latin1]{inputenc}
\usepackage{indentfirst}
\usepackage{color}
\usepackage{graphicx}
\usepackage{newlfont}
\usepackage{soul}
\usepackage{amssymb}
\usepackage{amsmath}
\usepackage{latexsym}
\usepackage{verbatim}
\usepackage{hyperref}
\usepackage{amsthm}
\usepackage{mathrsfs}
\usepackage{pifont}
\usepackage{textcomp}
\usepackage{tikz}
\usepackage{bm}

\theoremstyle{plain}

\newcommand{\IA}{\mathbb{A}}

\newcommand{\IC}{\mathbb{C}}

\newcommand{\IR}{\mathbb{R}}
\newcommand{\IZ}{\mathbb{Z}}

\newcommand{\IP}{\mathbb{P}}
\newcommand{\BO}{\mathcal{O}}

\newcommand{\BB}{\mathcal{B}}

\newcommand{\BX}{\mathcal{X}}
\newcommand{\BL}{\mathcal{L}}

\newcommand{\eeq}{\mathrel{\mathop:}=}

\newcommand{\Trop}{\operatorname{Trop}}

\newcommand{\Real}{\operatorname{Real}}

\newcommand{\V}{\operatorname{V}}
\newcommand{\Va}{\operatorname{V^{va}}}

\newcommand{\op}{\operatorname}
\newcommand{\T}{\operatorname{T}}
\newcommand{\RReal}{\operatorname{\mathbf{Real}}}

\newcommand{\mult}{\operatorname{mult}}
\newcommand{\pos}{\operatorname{pos}}
\newcommand{\I}{\operatorname{Ideal}}
\newcommand{\TAB}{\hspace*{10mm}}

\newtheorem{theorem}{Theorem}[section]

\newtheorem*{theorem*}{Theorem}
\newtheorem{proposition}[theorem]{Proposition}
\newtheorem*{proposition*}{Proposition}
\newtheorem*{lemma*}{Lemma}

\newtheorem{corollary}[theorem]{Corollary}
\newtheorem{lemma}[theorem]{Lemma}

\newtheorem{introteo}{Theorem}

\theoremstyle{definition}
\newtheorem{example}[theorem]{Example}
\theoremstyle{remark}
\newtheorem{remark}[theorem]{Remark}
\newtheorem{algoritmo}[theorem]{Algorithm}

\title[An Algorithm for the Realizability Problem for Families of Curves]{An Algorithm for the Tropical Realizability Problem for Families of Curves}
\author{Paolo Tripoli}
\address{School of Mathematical Sciences,
The University of Nottingham,
University Park,
Nottingham, NG7 2RD.}
\email{Paolo.Tripoli@nottingham.ac.uk}

\begin{document}

\begin{abstract}
Given a tropical fan curve $\Sigma$ and a family of algebraic curves $\BX \rightarrow \IA^k$ we define the realization locus $\Real_\Sigma \subseteq \IA^k$  as the set of fibers $X_a$ whose tropicalization is $\Sigma$. 
We produce an algorithm that describes the Zariski closure of $\Real_\Sigma$ by imposing algebraic conditions for each ray of $\Sigma$.
\end{abstract}

\maketitle

\section{Introduction}

The tropicalization map is a map from the set of embedded algebraic varieties to the set of weighted polyhedral complexes. It associates to a variety $X$ the tropical variety $\Trop(X)$. The tropical variety $\Trop(X)$ is called the tropicalization of $X$; it retains many properties of the original algebraic variety $X$.
Tropical geometry studies what algebro-geometric information of $X$ can be deduced from the polyhedral geometry and the combinatorics of $\Trop(X)$.
%
%
 One of the most elementary problems in tropical geometry is the realizability problem, which is the problem of deciding whether a given polyhedral complex $\Sigma$ can be obtained as the tropicalization of some algebraic variety.
 Some necessary conditions are well known (see~\cite[Section~2.5]{Speyer}): the polyhedral complex  $\Sigma$ must be rational, pure dimensional, balanced and connected in codimension one. These conditions are, however, far from being sufficient, and the realizability problem is very difficult in its full generality. In this paper we will describe an algorithm to solve a particular instance of the realization problem: we are given a tropical fan curve and a family of algebraic curves and we look for fibers of the family that tropicalize to the given tropical curve. 


Given an ideal $I\subseteq \IC[c_1, \ldots, c_k, x_1, \ldots, x_n]$ we consider the projection
\begin{equation}\label{eq:subsa}
	\begin{array}{ccc}
		\BX = \V(I)& \subseteq &\IA^k\times \IA^n \\
		\downarrow  & &\\
		\IA^k. & &
	\end{array}
\end{equation}
For fixed $a\in \IA^k$ we denote by $X_a$ the fiber of $\BX$ over $a$. We will refer to $\BX$ as a family of algebraic varieties over $\IA^k$, however no hypothesis of flatness is made on $\BX$.
The ideal $I_a$ of $X_a \subseteq \IA^n$ is the image of $I$ under the ring homomorphism
\begin{equation*}
	\begin{array}{cccc}
		i_a: & \IC[c_1, \ldots, c_k, x_1, \ldots, x_n] &\rightarrow & \IC[x_1, \ldots x_n] \\
				 & c_i																		&\mapsto 			& a_i\\
				 &  x_i																		 & \mapsto    & x_i.
	\end{array}
\end{equation*} 
%

We denote by $\T^n\subseteq \IA^n$ the torus $\T^n = \IA^n \setminus \V(x_1\cdots x_n)$ and, for an ideal $J\subseteq \IC[x_1, \ldots, x_n]$, we denote by $\Va(J)$ the vanishing locus $\V(J)\cap \T^n$ of $J$ in the torus. We study the tropicalization of the intersection $\Va(I_a)$ of the fibers of $\V(I)$ with the torus $\T^n$. We simply denote this tropicalization as $\Trop(I_a)$ or $\Trop(X_a)$. Both $\Trop(I_a)$ and $\Trop(X_a)$ strictly mean $\Trop(\Va(I_a))$. Moreover, we stress that we are working over the field of complex numbers, which is equipped with the trivial valuation. As a result, all tropical varieties we will consider are supported on a fan.

Let $\Sigma$ be a list of rational rays $\left(\rho_1, \ldots, \rho_l \right)$ with multiplicities $\left(m_1, \ldots, m_l\right)$. The reader is encouraged to think of $\Sigma$ as the (possibly partial) data of a tropical fan curve.
In this paper we introduce algorithms to describe the set
\begin{equation*} 
	\Real_{\Sigma} \eeq \{a\in \IA^k \mid  \dim (\Trop(I_a)) = 1 \text{ and } \mult(\rho_i, \Trop(I_a)) \geq m_i \text{ for }i=1,\ldots, l\}.  \\
\end{equation*}

Many instances of the realizability problem can be translated into this setting. For example, let $S = \V(f)$ be a surface in $\IA^3$, and let $\Sigma$ be a tropical curve contained in the tropicalization of $S$. Consider the relative realizability problem for $S$ and $\Sigma$, namely whether there exist an algebraic curve $C\subset S$ such that $\Sigma = \Trop(C)$. Suppose we further know that all the curves $C\subset S$ are complete intersection (this is true for very general $S$ as long as $\deg(S)\geq 4$, see~\cite{NL}). As a result, we are looking for a realization of $\Sigma$ of the form $C=\V(f, g)$. Moreover, as the degree of $C$ must equal the degree of $\Sigma$ as a tropical variety, the degree of $g$ is forced to be $d={\deg(f)}/{\deg(\Sigma)}$. Let $\V(I) \subset \IA^{\binom{3+d}{d}} \times \IA^3 $ be the family defined by the ideal $I=(f, \sum_{ijk} a_{ijk} x_1^ix_2^jx_3^k) \subset \IC[x_1, x_2, x_3, a_{ijk} \mid i+j+k \leq d]$.
The relative realizability problem for $S$ and $\Sigma$ is equivalent to the problem of realizability of $\Sigma$ over $\V(I)$.  A similar procedure, in theory, can be applied to all the cases of realizability problem in which it is known the number and the degrees of the generators of the ideal of the realization.

The main result of the paper is Algorithm~\ref{GLOBWEI}, which takes as input the ideal $I$ and the tropical curve $\Sigma$ and produces as output the ideal of the Zariski closure of $\Real_{\Sigma}$ in $\IA^k$. 
\begin{introteo}
Let $C$ be the output of Algorithm~\ref{GLOBWEI} with input $I$, $\Sigma$. Then $\V(C) = \overline{\Real(\Sigma)}$.
\end{introteo}
When $\Sigma$ is a tropical curve of degree $d$, and the intersection $\Va(I_a)$ of the fibers with the torus are known to be degree $d$ curves, the tropicalization $\Trop(I_a)$ of the fibers over $\Real_{\Sigma}$ equal $\Sigma$.

The algorithms described in this paper and a computation of all the examples have been implemented in a Macaulay2 (\cite{M2}) package available at \begin{center} 
\url{https://sites.google.com/view/paolotripoli/software/realizability}. \end{center}
The package has been developed as a proof of concept. Some of the routines presented in this paper are highly resource-demanding and their implementation in Macaulay2 is not fully optimized. 

The structure of the paper is the following. 
In Sections~\ref{sec2}~and~\ref{sec3} we restrict our attention to a single ray. In Section~\ref{sec2} we describe the set of fibers whose tropicalization contains the ray $\op{pos}(e_1)$. In Section~\ref{sec3} we describe the set of fibers whose tropicalization contains the ray $\op{pos}(e_1)$ with at least a fixed multiplicity $m$. 
Section~\ref{sec4} shows how to adapt these algorithms to describe, for any given ray $\rho$ and weight $m$, the set of fibers whose tropicalization contains the ray $\rho$ with multiplicity at least $m$. Finally in Section~\ref{sec5} we combine all the information coming from the rays of $\Sigma$ and describe Algorithm~\ref{GLOBWEI}.

\subsection*{Acknowledgments}  This work is part of the PhD thesis of the author. The author thanks his supervisor Diane Maclagan for close guidance and Christian Boehning and Alex Fink for useful comments.


\section{Local Set-Theoretical Computation}\label{sec2}

Let $I\subseteq \IC[c_1, \ldots, c_k, x_1, \ldots, x_n]$ be an ideal. In this section we study the set 
\begin{equation*}
	\op{Real_1} (I) \eeq \{ a\in \IA^k \mid \op{pos}(e_1)  \subseteq \Trop(I_a) \} 
\end{equation*}
of parameters $a$ such that the tropical variety $\Trop(I_a)$ contains $\op{pos}(e_1)$.

By Tevelev's Lemma (see \cite[Lemma~2.2]{Tev}) the tropical variety $\Trop(I_a)$ intersects the relative interior of the ray $\op{pos}(e_1)$ if and only if $\overline{\Va(I_a)}$ intersects the open torus orbit 
\begin{equation*}
	\BO_1 \eeq \{ (x_1, \ldots, x_n) \in \IA^n \mid x_1=0, \ x_i\neq 0 \text{ for } i=2,\ldots n \}.
\end{equation*}

We note that, as $\op{pos}(e_1)$ is a ray, $\Trop(I_a)$ intersects $\op{pos}(e_1)$ in its relative interior if and only if $\op{pos}(e_1)$ is completely contained in $\Trop(I_a)$. This gives the following equality of sets:
\begin{equation}\label{eq:real1}
\Real_1 (I)= \{ a \in \IA^k\ \mid \overline{\Va(I_a)}\cap \BO_1 \neq \emptyset \}.
\end{equation}
We also consider the following subset of $\IA^k\times \IA^n$:
\begin{equation}\label{eq:rreal1}
	\RReal_1 (I) \eeq \{(a,x)\in \IA^k\times \IA^n\ \mid \ x\in \overline{\Va(I_a)}\cap \BO_1\}.
\end{equation}

Whenever the ideal $I$ is clear from the context we will simply denote these sets as $\Real_1$ and $\RReal_1$. 
The sets $\Real_1$ and $\RReal_1$ are constructible sets and we denote by $\pi:\RReal_1 \rightarrow \Real_1$ the restriction of the projection $\IA^{k+n}\rightarrow \IA^k$.

\begin{figure}
\begin{center}
\begin{minipage}{0.3\textwidth}
	\begin{center}
		\begin{tikzpicture}[style=thick]
			\draw [fill=none] (-1,1.5cm) circle(0pt) node[above]{$\Sigma_1$};
			\draw [->](0,0cm) -- (1,1cm);
			\draw [->](0,0cm) -- (0,-1cm);
			\draw [->](0,0cm) -- (-1,0cm);
			\draw [fill=none](1,1cm) circle (0pt) node[right] {$(1,1)$};
			\draw [fill=none](0,-1cm) circle (0pt) node[left] {$(0,-1)$};
			\draw [fill=none](-1,0cm) circle (0pt) node[left] {$(-1,0)$}; 
			\draw [fill=white](0.5,0.5cm) circle (6pt) node[] {$1$};
			\draw [fill=white](-0.5,0cm) circle (6pt) node[] {$1$};
			\draw [fill=white](0,-0.5cm) circle (6pt) node[] {$1$};
		\end{tikzpicture}
	\end{center}
\end{minipage}
\begin{minipage}{0.3\textwidth}
	\begin{center}
		\begin{tikzpicture}[style=thick]
			\draw [fill=none] (-1,1.5cm) circle(0pt) node[above]{$\Sigma_2$};
			\draw [->](0,0cm) -- (0,1cm);
			\draw [->](0,0cm) -- (0,-1cm);
			\draw [fill=none](0,1cm) circle (0pt) node[right] {$(0,1)$};
			\draw [fill=none](0,-1cm) circle (0pt) node[left] {$(0,-1)$}; 
			\draw [fill=white](0,0cm) circle (6pt) node[] {$1$};
		\end{tikzpicture}
	\end{center}
\end{minipage}
\begin{minipage}{0.3\textwidth}
	\begin{center}
		\begin{tikzpicture}[style=thick]
			\draw [fill=none] (-1,1.5cm) circle(0pt) node[above]{$\Sigma_3$};
			\draw [->](0,0cm) -- (1,0cm);
			\draw [->](0,0cm) -- (-1,0cm);
			\draw [fill=none](1,0cm) circle (0pt) node[right] {$(1,0)$};
			\draw [fill=none](-1,0cm) circle (0pt) node[left] {$(-1,0)$}; 
			\draw [fill=none](0,-1cm) circle (0pt) node[left] {}; 
			\draw [fill=white](0,0cm) circle (6pt) node[] {$1$};
		\end{tikzpicture}
	\end{center}
\end{minipage}
\caption{The tropical curves of Example~\ref{eg:real1}.}
\label{fig:eg1}
\end{center}
\end{figure}
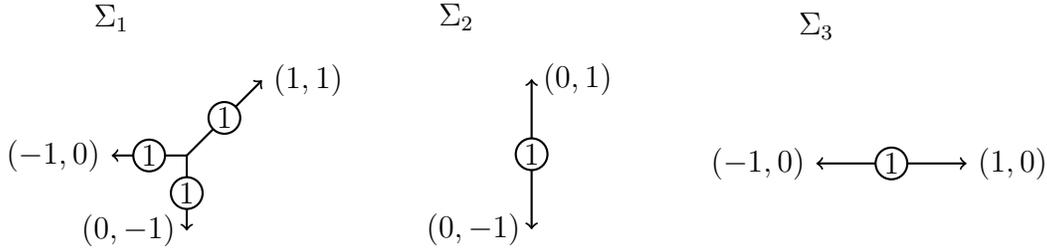

\begin{example}\label{eg:real1}
Consider the family $\V(I) \subseteq \IA^2 \times \IA^2$ where $I= (c_1x+c_2y+xy)\subseteq \IC[c_1, c_2, x, y]$. 
For $a=(a_1,a_2) \in \IC^2 \setminus \{(0,0)\}$, the tropical variety $\Trop(I_a)$ is depicted in Figure~\ref{fig:eg1}:
for $a_1,a_2\neq 0$ it is the curve $\Sigma_1$, for $a_1=0, a_2\neq 0$ it is the curve $\Sigma_2$, for $a_1\neq 0, a_2= 0$ it is the curve $\Sigma_3$. For $a=(0,0)$ the tropical variety $\Trop(I_a)$ is empty. In particular $e_1 \in \Trop(I_a)$ if and only if $a_1 \neq 0$ and $a_2=0$. Equivalently $\Real_1= \{(c_1, c_2) \mid c_1 \neq 0, c_2=0\}$.
\end{example}

The basic algebraic tools we will use are saturation and elimination theory; for a complete introduction see, for example, \cite[Section~14.1]{Eis} and \cite[Chapter~3, Section~15.10.6]{CLOS}. Given an ideal $I$ in a ring $R$, and an element $f\in R$, the saturation ideal $(I:f^\infty)$ is by definition 
\begin{equation*}
	(I:f^\infty) \eeq \{g\in R \mid f^ng \in I \text{ for some } n>0\}.
\end{equation*}
The ideal $I$ is said to be \emph{saturated} with respect to $f$ if $I = (I:f^\infty)$.

When $R$ is the polynomial ring $\IC[x_1, \ldots, x_n]$ we have $\V(I:f^\infty) = \overline{\V(I) \setminus \V(f)}$. In particular $\V(I: (x_1\cdots x_n)^\infty) = \overline{\Va(I)}$.

Given an ideal $I\subseteq \IC[c_1, \ldots, c_k, x_1, \ldots, x_n]$, the ideal $I \cap \IC[c_1, \ldots c_k]$ can be computed via elimination theory. The variety $\V(I \cap \IC[c_1, \ldots c_k])$ is the Zariski closure of the projection of $\V(I)\subseteq \IA^{k+n}$ to $\IA^k$.

Equations~\eqref{eq:real1} and \eqref{eq:rreal1} suggest a naive approach to compute the closure of $\RReal_1$ and $\Real_1$. As  $ \overline{\Va(I_a)}= V(I_a:(x_1\cdots x_n)^\infty)$ and $\BO_1 = (\V(x_1) \setminus \V(x_2\cdots x_n))$, the space $\RReal_1$ is given by points $(a,x)$ satisfying 
\begin{equation}\label{eq:realappr}
	x\in \overline{\Va(I_a)}\cap \BO_1 = V(I_a:(x_1\cdots x_n)^\infty) \cap (\V(x_1) \setminus \V(x_2\cdots x_n)).
\end{equation}

We now describe a naive algorithm to compute the ideals of $\RReal_1$ and $\Real_1$.
The naive algorithm approximates the ideal of $\RReal_1$ by replacing in Equation~\eqref{eq:realappr} the saturation of the ideal of the fiber $I_a$ with a saturation of the entire ideal $I$
\begin{equation}\label{eq:approx}
	(I_a:(x_1\cdots x_n)^\infty)  \leadsto (I:(x_1\cdots x_n)^\infty)_a,
\end{equation}
and, with the notation $A_2 =  (I_a:(x_1\cdots x_n)^\infty)+ (x_1)$, the set theoretic difference by $\V(x_2\cdots x_n)$ is replaced by a saturation
\begin{equation*}
	\V(A_2) \setminus \V(x_2\cdots x_n) \leadsto \V(A_2:(x_2\cdots x_n)^\infty).
\end{equation*}
Finally from the approximation $J=(A_2:(x_2\cdots x_n)^\infty)$ of $\I(\RReal_1)$, elimination theory is used to approximate the ideal $\I(\Real_1)$
\begin{equation*}
	\I(\Real_1) \leadsto J\cap \IC[c_1, \ldots, c_k].
\end{equation*}
To sum up, the naive algorithm performs the following operations:
\begin{enumerate}
	\item $A_1=( I:x_1^\infty )$
	\item $A_2= A_1+ (x_1)$
	\item $J=(A_2:(x_2\cdots x_n)^\infty)$
	\item $C=J\cap \IC[c_1,\ldots, c_k]$
	\item \texttt{return} $(C, J)$.
\end{enumerate}


The naive approach does not work well as the following example shows.

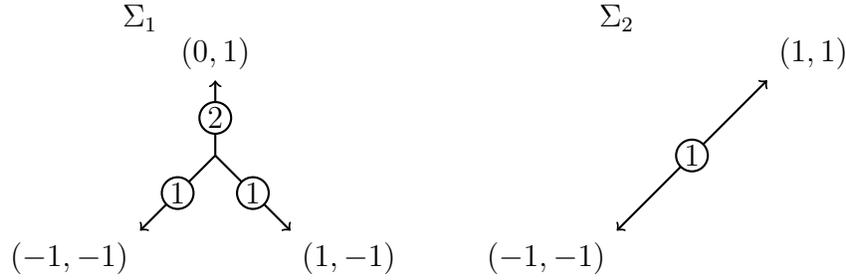
\begin{figure}
\begin{center}
\begin{minipage}{0.3\textwidth}
	\begin{center}
		\begin{tikzpicture}[style=thick]
			\draw [fill=none] (-1,1.5cm) circle(0pt) node[above]{$\Sigma_1$};
			\draw [->](0,0cm) -- (0,1cm);
			\draw [->](0,0cm) -- (-1,-1cm);
			\draw [->](0,0cm) -- (1,-1cm);
			\draw [fill=none](0,1cm) circle (0pt) node[above] {$(0,1)$};
			\draw [fill=none](-1,-1cm) circle (0pt) node[below left] {$(-1,-1)$};
			\draw [fill=none](1,-1cm) circle (0pt) node[below right] {$(1,-1)$}; 
			\draw [fill=white](0,0.5cm) circle (6pt) node[] {$2$};
			\draw [fill=white](0.5,-0.5cm) circle (6pt) node[] {$1$};
			\draw [fill=white](-0.5,-0.5cm) circle (6pt) node[] {$1$};
		\end{tikzpicture}
	\end{center}
\end{minipage}
\hspace{0.1\textwidth}
\begin{minipage}{0.3\textwidth}
	\begin{center}
		\begin{tikzpicture}[style=thick]
			\draw [fill=none] (-1,1.5cm) circle(0pt) node[above]{$\Sigma_2$};
			\draw [->](0,0cm) -- (1,1cm);
			\draw [->](0,0cm) -- (-1,-1cm);
			\draw [fill=none](1,1cm) circle (0pt) node[above right] {$(1,1)$};
			\draw [fill=none](-1,-1cm) circle (0pt) node[below left] {$(-1,-1)$}; 
			\draw [fill=white](0,0cm) circle (6pt) node[] {$1$};
		\end{tikzpicture}
	\end{center}
\end{minipage}
\caption{The tropical curves of Example~\ref{eg:naive}.}
\label{fig:eg2}
\end{center}
\end{figure}

\begin{example}\label{eg:naive}
Let us consider the ideal $I=(c_2x^2+c_2xy+c_1)\subseteq \IC[c_1,c_2,x,y]$. 
The tropical variety $\Trop(I_{(a_1,a_2)})$ for $a_2 \neq 0$ is depicted in Figure~\ref{fig:eg2}: it is the tropical curve $\Sigma_1$ for  $a_1,a_2 \neq 0$ and the tropical curve $\Sigma_2$ for $a_1=0, a_2\neq 0$. The tropical variety $\Trop(I_{(a_1,a_2)})$ is empty  for $a_1\neq 0$ and $a_2=0$, and it is the whole $\IR^2$ when $a_1=a_2=0$.
In particular we have $\Real_1 = \{(0,0)\}$. We would wish the naive Algorithm to produce as output the ideal $(c_1, c_2)$.
However, running it we get
\begin{itemize}
	\item $A_1=(I:x^\infty)=I$,
	\item $A_2=A_1+(x) = I+(x) = (c_1, x)$,
	\item $J= (A_2 : y^\infty ) = (c_1, x)$,
	\item $C= J \cap \IC[c_1, c_2] = (c_1)$. 
\end{itemize}
\end{example}

This happens because performing the saturation on the entire family as we do in Equation~\eqref{eq:approx} is a weaker operation than performing it on the fiber. We will give a precise statement of this fact in Corollary~\ref{satulemma0}.

We note that, in Example~\ref{eg:naive}, while the ideal $I=(c_2x^2+c_2xy+c_1)$ is saturated with respect to $(x)$, the ideal $I+(c_1)= (c_2x^2+c_2xy, c_1)$ is not. Actually, we can observe that running again the naive Algorithm on $I'=I + (c_1)$, we get as output $C'= (1)$ as desired. This is a general fact: the ideal of $\Real_1$ can be computed by a repeated application of the naive Algorithm. More precisely, we introduce the following algorithm.

\begin {algoritmo}\label{LOCALG}
\leavevmode \\
{\bf {Input:}}  
			\begin{itemize}
				\item[$I$,] ideal in $\IC[c_1,\ldots, c_k, x_1, \ldots, x_n]$. 
			\end{itemize}
{\bf {Output:}}  
			\begin{itemize}
				\item[$C$,] ideal in $\IC[c_1, \ldots, c_k]$ satisfying $\V(C)=\overline{\Real_1}(I)$,
				\item[$J$,] ideal in $\IC[c_1, \ldots, c_k, x_1, \ldots x_n]$ satisfying $\V(J_a)= \overline{(\RReal_1(I))_a}$ for $a$ general in $\V(C)$.
			\end{itemize}
\{ \\
$C=(0) \subseteq \IC[c_1, \ldots, c_k]$; \\
$J_{\text{old}}=(1)$; \\
$J=(0) \subseteq \IC[c_1, \ldots, c_k, x_1, \ldots, x_n]$; \\
\texttt{while} ($J_{\text{old}} \neq J$) \texttt{do} $\{$ \\
\TAB $J_{\text{old}}=J$; \\
\TAB $A_1=( ( I + C):x_1^\infty)$; \\
\TAB $A_2= A_1+(x_1)$; \\
\TAB $J=(A_2:(x_2\cdot\ldots\cdot x_n)^\infty)$; \\
\TAB $C=\op{Radical}(J\cap \IC[c_1,\ldots, c_k])$; \\
\TAB $\}$ \\
\texttt{return} $(C, J)$; \\
\}
\proof[Proof of correctness]
This follows from Lemma~\ref{lem:localgrest}, Proposition~\ref{prop:torusfibers} and Theorem~\ref{localgteo2} below.
\endproof
\end{algoritmo}

For the proofs of the following results we denote by $A_1^{(i)},\ A_2^{(i)},\ J^{(i)},\ C^{(i)}$ the ideals $A_1,\ A_2,\ J,\ C$ as computed at the $i$\textsuperscript{th} iteration of the \texttt{while} loop of Algorithm \ref{LOCALG}. 

\begin{lemma}\label{lem:localgrest}
For every $i>0$ we have $C^{(i)}\subseteq C^{(i+1)}$ and $J^{(i)} \subseteq J^{(i+1)}$.
In particular, by Noetherianity, Algorithm \ref{LOCALG} terminates after a finite number of steps.
\proof
Fix $i>0$  and consider the ideal $C^{(i)}\subseteq \IC[c_1, \ldots, c_k]$. We have $ C^{(i)} \subseteq A_1^{(i+1)} \subseteq A_2^{(i+1)} \subseteq J^{(i+1)}$, and hence $C^{(i)}  \subseteq C^{(i+1)} =  \op{Radical}(J^{(i+1)} \cap \IC[c_1, \ldots, c_k])$. 

From $C^{(i)}\subseteq C_1^{(i+1)}$ it follows that $A_1^{(i)}\subseteq A_1^{(i+1)}$, $A_2^{(i)}\subseteq A_2^{(i+1)}$ and, finally, $J^{(i)}\subseteq J^{(i+1)}$.
\endproof
\end{lemma}

\begin{remark}\label{rm:localgtriv}
We have $I\subseteq A^{(0)}_1 \subseteq A^{(0)}_2 \subseteq J_1^{(0)}$. In particular, using Lemma~\ref{lem:localgrest} we have that $\V(J)\subseteq \V(I)$.
\end{remark}

We recall the following fact about saturations.

\begin{lemma}\label{satulemma0.5}
Let $J,K\subseteq \IC[y_1, \ldots, y_m]$ be ideals and $f\in\IC[y_1,\ldots, y_m]$. Then $(J:f^\infty) + K \subseteq (J+K:f^\infty)$ and $((J:f^\infty)+ K:f^\infty) = (J+K:f^\infty)$.
\end{lemma}

We will mostly use Lemma~\ref{satulemma0.5} in the form of the following corollary.

\begin{corollary}\label{satulemma0}
Let $J\subseteq \IC[c_1,\ldots, c_k, x_1,\ldots, x_n]$ be an ideal and $f\in\IC[x_1,\ldots, x_n]$, then, for every $a\in \IC^k$, $(J:f^\infty)_a\subseteq (J_a:f^\infty)$ and $(J_a:f^\infty)=((J:f^\infty)_a:f^\infty)$.
\proof
This follows from Lemma~\ref{satulemma0.5}, after identifying the ideal $J_a \subseteq \IC[x_1, \ldots, x_n]$ with the ideal $J + (c_1-a_1, \ldots, c_k-a_k) \subseteq \IC[c_1,\ldots, c_k, x_1,\ldots, x_n]$.
\endproof
\end{corollary}

We are now ready to prove the first relation between the ideals $C$ and $J$ and the spaces $\Real_1$ and $\RReal_1$, namely that $\Real_1 \subseteq \V(C)$ and $\RReal_1 \subseteq \V(J)$.

\begin{proposition}\label{localgprop1}
Let $(C, J)$ be the output of Algorithm \ref{LOCALG} with input $I$. Then we have $\Real_1 \subseteq \V(C)$ and $\RReal_1 \subseteq \V(J)$.
\proof
As $\Real_1 = \pi(\RReal_1)$ and $\pi(\V(J)) \subseteq \V(C)$, it is enough to prove that $\RReal_1 \subseteq \V(J)$.
We prove by induction that $\RReal_1 \subseteq \V(J^{(i)})$ for every $i\geq 0$, the base case $i=0$ being trivial. 
Fix $(a,x)\in \RReal_1$ and assume $(a,x) \in \V(J^{(i)})$. 
From $(a,x) \in \V(J^{(i)})$ it follows that $a\in \V(C^{(i)})$. 
Fix a polynomial $f \in C^{(i)}$. If we regard $f$ as a polynomial in $\IC[c_1, \ldots, c_k]$ we have $f(a)=0$ as $a \in \V(C^{(i)})$. Equivalently, if we regard $f$ as a polynomial in $\IC[c_1, \ldots, c_k, x_1, \ldots, x_n]$, we have that $f_a = i_a(f)$ is the zero polynomial in $\IC[x_1, \ldots, x_n]$, where $i_a$ is defined by Equation~\ref{eq:subsa}.
As a result, for every $f\in C^{(i)}$, we have $f_a = 0$, which means that $C^{(i)}_a = (0)$. 
It follows that $I_a = I_a + C^{(i)}_a = (I + C^{(i)})_a$. 
As we assumed that $(a,x)\in \RReal_1$ we get, by the definition of $\RReal_1$, that $x$ is in $\overline{\Va (I_a)}$. We have $\overline{\Va(I_a)} = \overline{\Va (I_a:x_1^\infty)}= \overline{\Va ((I+ C^{(i)})_a:x_1^\infty)}$. 
In particular, $x\in \V ((I + C^{(i)})_a:x_1^\infty)$. 
Corollary~\ref{satulemma0} implies that $x \in \V(((I + C^{(i)}):x_1^\infty)_a) = \V((A_1^{(i+1)})_a)$. 
By the definition of $\RReal_1$ we have $x\in \BO_1 = \V(x_1) \setminus \V(x_2\cdots x_n)$.
From $x\in \V(x_1)$ it follows that $x \in \V ((A_1^{(i+1)})_a + (x_1)) = \V((A_2^{(i+1)})_a)$.
From $x\notin \V(x_2\cdots x_n)$ it follows that $x \in \V ((A_2^{(i+1)})_a) \setminus \V(x_2\cdots x_n) \subseteq \V ((A_2^{(i+1)})_a:(x_2\cdots x_n)^\infty)$. 
Again by Corollary~\ref{satulemma0}, we have that $x\in \V ((A_2^{(i+1)}:(x_2\cdots x_n)^\infty)_a)= \V(J^{(i+1)}_a)$. 
This last condition is equivalent to $(a,x)\in \V(J^{(i+1)})$, which concludes the proof.
\endproof
\end{proposition}

Corollary~\ref{satulemma0} shows that the saturation ideal of a fiber is a contained in the fiber of the saturation of the entire family. The following lemma shows that this containment is generally an equality.

\begin{lemma}\label{satulemma1}
Let $J\subseteq R \eeq \IC[c_1,\ldots, c_k, x_1,\ldots, x_n]$ be an ideal and suppose that $C=J\cap \IC[c_1, \ldots, c_k]$ is radical. Suppose that $J$ is saturated with respect to $x_1$. Then there exists an open dense subset $A\subseteq \V(C)$ such that, for every $a \in A$,  $J_a$ is saturated with respect to $x_1$.
\end{lemma}

The proof of Lemma~\ref{satulemma1} is postponed to Section~\ref{Ssec:HiDim}, where we describe an algorithm to compute the complement of the set $A$.

\begin{corollary}\label{satucoro2}
Let $J\subseteq \IC[c_1,\ldots, c_k, x_1,\ldots, x_n]$ be an ideal and suppose that $C=J\cap \IC[c_1, \ldots, c_k]$ is radical. Then there exists an open dense subset $A\subseteq \V(C)$ such that, for every $a \in A$,  $(J:x_1^\infty)_a=(J_a:x_1^\infty)$.
\proof
Applying Lemma~\ref{satulemma1} to $(J:x_1^\infty)$ we get that there exists an open dense subset $A\subseteq \V(C)$ such that, for every $a \in A$,  $(J:x_1^\infty)_a$ is saturated with respect to $x_1$. Equivalently  $(J:x_1^\infty)_a =( (J:x_1^\infty)_a:x_1^\infty)$ and the result now follows from Corollary~\ref{satulemma0}.
\endproof
\end{corollary}

Corollary~\ref{satucoro2} allows us to show that the fibers of $\V(J)$ agree with the fibers of $\RReal_1$. This is the content of the following proposition.

\begin{proposition}\label{prop:torusfibers}
Let $(C, J)$ be the output of Algorithm \ref{LOCALG} with input $I$. Then for $a \in \V(C)$ general we have $\V(J_a) = \overline{(\RReal_1)_a}$.
\proof
We have 
\begin{align*}
	(\RReal_1)_a& = \overline{\Va(I_a)} \cap \BO_1\\
						&=\V \left(I_a: (x_1\cdots x_n)^\infty\right) \cap \left( \V(x_1)  \setminus \V(x_2\cdots x_n) \right) \\
						&=\V \left(\left(I_a: (x_1\cdots x_n)^\infty\right) + (x_1\right))  \setminus \V(x_2\cdots x_n)  \\					
						&=\V \left(\left(\left(I_a: (x_1\cdots x_n)^\infty\right)+ (x_1\right)): (x_2\cdots x_n)^\infty \right)  \setminus \V(x_2\cdots x_n)  \\
						&=\V \left(\left(\left(I_a: x_1^\infty\right)+ (x_1\right)): (x_2\cdots x_n)^\infty \right)  \setminus \V(x_2\cdots x_n),
\end{align*}
where the last equality comes from Lemma~\ref{satulemma0.5}. It follows that 
\begin{equation*}
	\overline{(\RReal_1)_a} = \V \left(\left(\left(I_a: x_1^\infty\right)+( x_1\right)): (x_2\cdots x_n)^\infty \right).
\end{equation*}
By construction we have
\begin{equation}\label{Eq.CinJ}
	J =  \left(\left(\left(I+C: x_1^\infty\right)+ (x_1\right)): (x_2\cdots x_n)^\infty \right)
\end{equation}
and, therefore,
\begin{equation*}
	J_a =  \left(\left(\left(I+C: x_1^\infty\right)+ (x_1\right)): (x_2\cdots x_n)^\infty \right)_a.
\end{equation*}
We now show that the ideal $J$ satisfies the assumptions of Corollary~\ref{satucoro2}, namely that $C'= J \cap \IC[c_1, \ldots, c_k]$ is radical. By construction we have $C = \op{Radical}(C')$. Moreover, by Equation~\eqref{Eq.CinJ}, we have $C\subseteq J$ and therefore $C \subseteq C'$. This shows that $C'= C$ and in particular $C' = J \cap \IC[c_1, \ldots, c_k]$ is radical. Similarly we show that $(I+C) \cap \IC[c_1,\ldots c_k] =C$, so that $I+C$ satisfies the assumptions of Corollary~\ref{satucoro2}.
We $C \subseteq (I+C) \cap \IC[c_1,\ldots c_k]$ because  $C \subseteq I+C$ and, since $I\subseteq J$, we also have $I \cap \IC[c_1, \ldots, c_k] \subseteq J \cap \IC[c_1, \ldots, c_k] = C$. This shows that $C=I \cap \IC[c_1, \ldots, c_k]$.

By applying Corollary~\ref{satucoro2} to $J$ and to $I+C$ we get that, for $a\in \V(C)$ general,
\begin{align*}
	J_a =  &\left(\left(I+C: x_1^\infty\right)+ (x_1)\right)_a: (x_1\cdots x_n)^\infty \\
		&=\left(\left(I+C: x_1^\infty\right)_a+ (x_1)\right): (x_1\cdots x_n)^\infty \\
		&=\left(\left((I+C)_a: x_1^\infty\right)+ (x_1)\right): (x_1\cdots x_n)^\infty.\\
\end{align*}
where the first equality is Corollary~\ref{satucoro2} applied to $J$, and the third equality is Corollary~\ref{satucoro2} applied to $I+C$.
The result  follows from $I_a = (I+C)_a$.
\endproof
\end{proposition}

We are now ready to prove the main result of this section.

\begin{theorem}\label{localgteo2}
Let $(C, J)$ be the output of Algorithm \ref{LOCALG} with input $I$. Then we have $\overline{\Real}_1 = \V(C)$. 
\proof
Let $A\subseteq \V(C)$ be an open dense subset on which the condition of Proposition~\ref{prop:torusfibers} holds. For $a\in A$ we have that $\V(J)_a$ is empty if and only if $(\RReal_1)_a$ is empty. It follows that $\pi(\V(J))\cap A = \pi(\RReal_1) \cap A$. In particular $\Real_1 = \pi(\RReal_1) \supset \pi(\V(J)) \cap A$ which is dense in $\V(C)$. This shows  $\overline{\Real}_1 = \V(C)$ .
%
\endproof
\end{theorem}

											\section{Local Weighted Computation}\label{sec3}

The aim of this section is to describe the Zariski closure of the set
\begin{align*}
\Real_{1,m} (I) \eeq \{a\in \IA^k \mid  \dim (\Trop(I_a)) = 1 \text{ and } \mult(\op{pos}(e_1), \Trop(I_a)) \geq m\}.
\end{align*}
We will often omit the ideal $I$ and simply write $\Real_{1,m}$.

\begin{figure}
\begin{center}
\begin{minipage}{0.3\textwidth}
	\begin{center}
		\begin{tikzpicture}[style=thick]
			\draw [fill=none] (-1,1.5cm) circle(0pt) node[above]{$\Sigma_1$};
			\draw [->](0,0cm) -- (0,1cm);
			\draw [->](0,0cm) -- (-1,-1cm);
			\draw [->](0,0cm) -- (1,0cm);
			\draw [fill=none](0,1cm) circle (0pt) node[above] {$(0,1)$};
			\draw [fill=none](-1,-1cm) circle (0pt) node[below left] {$(-1,-1)$};
			\draw [fill=none](1,0cm) circle (0pt) node[right] {$(1,0)$};
			\draw [fill=white](0.5,0cm) circle (6pt) node[] {$1$};
			\draw [fill=white](0,0.5cm) circle (6pt) node[] {$1$};
			\draw [fill=white](-0.5,-0.5cm) circle (6pt) node[] {$1$};
		\end{tikzpicture}
	\end{center}
\end{minipage}
\begin{minipage}{0.3\textwidth}
	\begin{center}
		\begin{tikzpicture}[style=thick]
			\draw [fill=none] (-1,1.5cm) circle(0pt) node[above]{$\Sigma_2$};
			\draw [->](0,0cm) -- (0,1cm);
			\draw [->](0,0cm) -- (-1,-1cm);
			\draw [->](0,0cm) -- (1,0cm);
			\draw [fill=none](0,1cm) circle (0pt) node[above] {$(0,1)$};
			\draw [fill=none](-1,-1cm) circle (0pt) node[below left] {$(-1,-1)$};
			\draw [fill=none](1,0cm) circle (0pt) node[right] {$(1,0)$};
			\draw [fill=white](0.5,0cm) circle (6pt) node[] {$2$};
			\draw [fill=white](0,0.5cm) circle (6pt) node[] {$2$};
			\draw [fill=white](-0.5,-0.5cm) circle (6pt) node[] {$2$};
		\end{tikzpicture}
	\end{center}
\end{minipage}
\begin{minipage}{0.3\textwidth}
	\begin{center}
		\begin{tikzpicture}[style=thick]
			\draw [fill=none] (-1,1.5cm) circle(0pt) node[above]{$\Sigma_3$};
			\draw [->](0,0cm) -- (0,1cm);
			\draw [->](0,0cm) -- (-1,-1cm);
			\draw [->](0,0cm) -- (1,0cm);
			\draw [->](0,0cm) -- (0,-1cm);
			\draw [fill=none](0,1cm) circle (0pt) node[above] {$(0,1)$};
			\draw [fill=none](-1,-1cm) circle (0pt) node[below left] {$(-1,-1)$};
			\draw [fill=none](1,0cm) circle (0pt) node[right] {$(1,0)$};
			\draw [fill=none](0,-1cm) circle (0pt) node[below] {$(0,-1)$};
			\draw [fill=white](0.5,0cm) circle (6pt) node[] {$1$};
			\draw [fill=white](0,0.5cm) circle (6pt) node[] {$2$};
			\draw [fill=white](-0.5,-0.5cm) circle (6pt) node[] {$1$};
			\draw [fill=white](0,-0.5cm) circle (6pt) node[] {$1$};
		\end{tikzpicture}
	\end{center}
\end{minipage}
\caption{The tropical varieties of Example~\ref{eg:6}.}
\label{fig:eg6}
\end{center}
\end{figure}
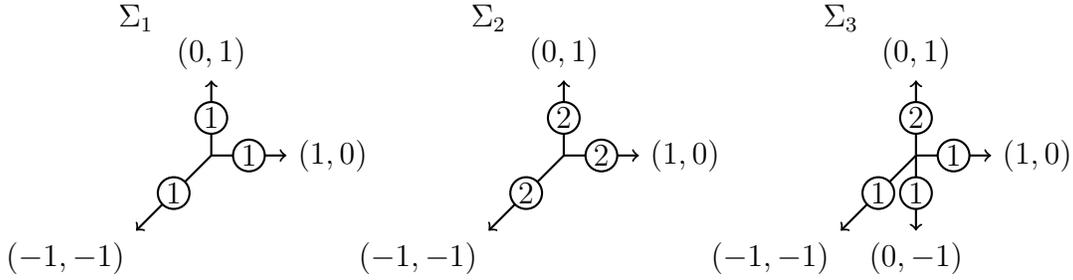

\begin{example}\label{eg:6}
Consider in $\IC[c_1,c_2,x,y]$ the following ideals: $I_1=(c_1, x+y+1)$, $I_2=(c_2, x+7y+4)$, $I_3=(c_1-c_2-1, x+c_1y+1)$. Let $I\subseteq \IC[c_1,c_2,x,y]$ be the intersection $I = I_1 \cap I_2 \cap I_3$. The tropicalization of the fiber over $a$ of the projection of $\V(I)$ to $\IC^2$ is empty outside $\V(c_1c_2(c_1-c_2-1)) \subseteq \IC^2$ and it is depicted in Figure~\ref{fig:eg6} for $a \in \V(c_1c_2(c_1-c_2-1))$: it is the tropical curve $\Sigma_1$ for $a$ in $\V(c_1c_2(c_1-c_2-1))\setminus \{ (0,0),(0,-1), (1,0)\}$, it is the tropical curve $\Sigma_2$ for $a$ in $\{(0,0), (1,0)\}$ and it is the tropical curve $\Sigma_3$ for $a=(0,-1)$. In particular the realization space $\Real_{1, 2}$ consists of the points $(0,0)$ and $(1,0)$. 
\end{example}

The computation of the multiplicity of the ray $\op{pos}(e_1)$ in the tropicalization of a very affine curve $X\subseteq \T^n$ can be based on the following generalization of Tevelev's Lemma (see, for example, \cite[Section~2.5]{Speyer}).

\begin{theorem}\label{prop:keymult}
Let $X\subseteq \T^n$ be an irreducible curve. Then $\mult(\op{pos}(e_1), \Trop(I))$ equals the cardinality, counted with multiplicity, of $\overline{X} \cap \BO_1$.
\end{theorem}

To compute the Zariski closure of $\Real_{1,m}$ we first use Algorithm~\ref{LOCALG} to approximate $\Real_1$ and $\RReal_1$. 
Since we do not assume any flatness condition on the ideal $I$, the next step is to compute the set of parameters $a$ for which $\Trop(I_a)$ has dimension $1$. This is carried out in Section~\ref{Ssec:HiDim}.
 In Section~\ref{Ssec:HiHil} we use Theorem~\ref{prop:keymult} to study the multiplicity of $\op{pos}(e_1)$ in the tropicalization of the fibers $\V(I_a)$.
 Finally, in Section~\ref{Ssec:HiGlu}, we combine all this information to describe $\Real_{1,m}$.


	\subsection{Higher dimensional fibers}\label{Ssec:HiDim}

In this section we describe the set of parameters $a$ such that the tropicalization of the fiber $\Trop(I_a)$ has dimension at most $1$. The dimension of $\Trop(I_a) = \Trop(\Va(I_a))$ equals the dimension of $\Va(I_a)$. Consider an irreducible component $Y$ of $\V(I)$ and let $Z$ be the closure of the projection of $Y$ to $\IA^k$. Let $d$ be the dimension of $Z$ and let $d+l$ be the dimension of $Y$. By~\cite[Theorem~2 of Chapter~1 Section~8]{Mumred}, for every $a\in Z$, every component of the fiber $Y_a$ has dimension at least $l$. Now suppose that $l\geq 2$ and fix a point $a \in Z$. If the fiber $Y_a$ is not empty its dimension, and thus the dimension of $\V(I_a)$, is at least $2$. Denote by $\overline{Y}$ the closure of $Y$ in $\IP^n\times \IA^k$. By the Main Theorem of Elimination Theory (see~\cite[Theorem~14.1]{Eis}), the projection $\pi(\overline{Y})$ of $\overline{Y}$ is closed in $\IA^k$. As a result, $\pi(\overline{Y})$ contains $Z$ and, in particular, the fiber $(\overline{Y})_a$ is not empty. As a result, a necessary condition for $\dim(\Va(I_a)) \leq 1$ is that the fiber $(\overline{Y})_a \subseteq \IP^n$ does not intersect the torus.
Finally, a necessary condition for the fiber $(\overline{Y})_a$ not to intersect the torus is that the ideal of $(\overline{Y})_a$ is not saturated.


\begin{example}
Consider the ideals  $I_1 = (x+cy+cz+c) \subseteq \IC[c,x,y,z]$, $I_2 = (x+2y-z+3,y+2z+2)\subseteq \IC[c, x,y,z]$, and let $I$ be the intersection $I_1\cap I_2$. The tropicalization of the fiber $\V(I_a)$, for $a\neq 0$, is the standard tropical plane. However, for $a=0$, the tropicalization of $\Trop(\V(I_0))$ is the tropical standard line, because $\V(x)$ does not intersect the torus.
\end{example}

Let $J\subseteq \IC[c_1,\ldots, c_k, x_1,\ldots, x_n]$ be an ideal saturated with respect to $x_n$ and suppose that $C=J\cap \IC[c_1, \ldots, c_k]$ is prime. 
A \emph{comprehensive} Gr\"obner basis is a Gr\"obner basis $\BL$ of $J$ with the extra property that $\BL_a = \{f_a \mid a\in \BL\}$ is a Gr\"obner basis of $J_a$ for every $a \in \IA^k$. Comprehensive Gr\"obner bases were first introduced, together with an algorithm to compute them, in \cite{ComprehensiveGB}. A more efficient algorithm to compute comprehensive Gr\"obner bases is described in \cite{CompGBAlgo}.
Let $\BL$ be a comprehensive Gr\"obner basis of $J$ with respect to the graded reverse lexicographical order. For each $f\in \BL$ write $f = \sum_{i\geq 0} f_i x_n^i$, where $f_i \in \IC[c_1, \ldots, c_k, x_1, \ldots, x_{n-1}]$ for every $i$. We can assume that $f_0 \notin J$ or $f=f_0$. Indeed, if $f = f_0 + x_n f'$ with $f_0 \in J$, then $x_nf' \in J$ and hence, by saturation, $f'\in J$. Therefore we can obtain a new comprehensive Gr\"obner basis of $J$ by removing $f$ from $\BL$ and adding to it $f_0$ and $f'$. Let $\BL_1$ be the subset of $\BL$ that consists of polynomials $f$ that satisfy $f \neq f_0$. 
For every $f\in \BL_1$ we denote by $M_f \subseteq \IC[c_1,\ldots, c_k]$ the ideal generated by the coefficients of $f_0$, where $f_0$ is now regarded as a polynomial in $(\IC[c_1, \ldots, c_k])[x_1, \ldots, x_{n-1}]$. Let $L\subseteq \IC[c_1,\ldots, c_k]$ be the ideal $\bigcap_{f\in \BL_1} M_f$.

\begin{proposition}\label{Prop:satulemmahard}
Let $J\subseteq \IC[c_1,\ldots, c_k, x_0,\ldots, x_n]$ be an ideal homogeneous with respect to the $x_i$'s and saturated with respect to $x_n$. Suppose that $C=J\cap \IC[c_1, \ldots, c_k]$ is prime. Then, with the notation introduced above, for every $a \in \V(C) \setminus \V(L+C)$, the ideal $J_a$ is saturated with respect to $x_n$, moreover the set $\V(C) \setminus \V(L+C)$ is dense in $\V(C)$.
\proof
Let $\BL$ be the comprehensive Gr\"obner basis of $J$ introduced above and fix $a\in \V(C)$. We have that $\BL_a \eeq \{f_a \mid f\in\BL\}$ is a Gr\"obner basis of $J_a$. Let $\BL'_a$ be the set obtained by dividing every element $f_a \in \BL_a$ by the highest power of $x_n$ that divides $f_a$. By \cite[Lemma~12.1]{BasesandPolytopes} $\BL'_a$ is a Gr\"obner basis of $(J_a:x_n^\infty)$. Fix $f\in \BL$. If $f\notin \BL_1$ then $f\in \IC[c_1, \ldots, c_k, x_0, \ldots, x_{n-1}]$ and $x_n$ does not divide $f_a$ for every $a$. If $f\in \BL_1$ and $a\notin \V(M_f)$, by construction, $f_a$ is not divisible by $x_n$. To sum up, for $a \in \V(C) \setminus \V(L+C)$ none of the polynomials in $\BL_a$ is divisible by $x_n$. As a result we have $\BL_a = \BL'_a$ and therefore $J_a = (J_a:x_n^\infty)$, since they are generated by the same set of polynomials.   
To conclude we show that the complement of $\V(C+M_f)$ is dense in $\V(C)$ for every $f\in \BL_1$. Let $f\in \BL_1$ and write $f_0 = \sum g_u x^u$ with $g_u \in \IC[c_1, \ldots, c_k]$ for every $u$. By the definition of $\BL_1$, we have $f\neq f_0$. By our assumption on $\BL$ this implies that $f_0 \notin J$ and therefore not all the $g_u$ can be in $C$. As $C$ is prime, it follows that $\V(C+M_f)$ is a proper closed set in $\V(C)$.
\endproof
\end{proposition}

\proof[Proof of Lemma~\ref{satulemma1}]
This is a weak formulation of Proposition~\ref{Prop:satulemmahard}.
\endproof

Proposition~\ref{Prop:satulemmahard} allows us to write the following algorithm.

\begin {algoritmo}\label{satualgo}
\leavevmode \\
{\bf {Input:}}  
			\begin{itemize}
				\item[$J$,] ideal in $\IC[c_1,\ldots, c_k, x_0, \ldots, x_n]$ homogeneous with respect to the $x_i$'s and saturated with respect to $x_n$.
			\end{itemize}
{\bf {Output:}}  
			\begin{itemize}
				\item[$C+L$,] ideal in $\IC[c_1, \ldots, c_k]$ satisfying $\V(C+L) = \{a \mid J_a \neq (J_a:x_n)\}$.
			\end{itemize}
\{ \\
$\op{TerminationCondition}= \op{False}$;\\
\texttt{while} ($\op{TerminationCondition}= \op{False}$) \texttt{do} \\
\TAB $\{$\\
\TAB $C= J \cap \IC[c_1, \ldots, c_k]$; \\
\TAB $C= \op{Radical}(C)$; \\
\TAB $J= J +C$;\\
\TAB compute a comprehensive Gr\"obner basis $\BB$ of $J$ with respect to the graded reverse lexicographical order; \\
\TAB $\BL=\{\}$;\\
\TAB \texttt{while}($\BB \neq \emptyset$) \texttt{do}\\
\TAB \TAB \{ \\
\TAB \TAB let $f$ be the first element of $\BB$;\\
\TAB \TAB write $f = f_0+ x_nf'$, where $f_0 \in \IC[c_1, \ldots, c_k, x_0, \ldots, x_{n-1}]$;\\
\TAB \TAB \texttt{if} ($f_0 \in J$ \texttt{AND} $f'\neq 0$) \texttt{then} \\
\TAB \TAB \TAB \{add $f_0$ to $\BL$; remove $f$ from $\BB$; add $f'$ to $\BB$\}; \\
\TAB \TAB \texttt{else}\\
\TAB \TAB \TAB  \{add $f$ to $\BL$; remove $f$ from $\BB$\}; \\ 
\TAB \TAB \} \\
\TAB $L = C$; \\
\TAB $\BL_1 = \{ f \in \BL \mid f \notin \IC[c_1, \ldots, c_k, x_0, \ldots, x_{n-1}]\}$; \\
\TAB \texttt{for} $f$ \texttt{in} $\BL_1$ \texttt{do} \\
\TAB \TAB  $\{$ \\
\TAB \TAB  write $f$ as $\sum a_i x_n^i$ with $a_i \in \IC[c_1, \ldots, c_k] [x_1, \ldots, x_{n-1}]$; \\
\TAB \TAB 		let $M_f \subseteq \IC[c_1, \ldots, c_k]$ be the ideal generated by the coefficients of $a_0$; \\
\TAB \TAB 		$L = L\cap M_f$; \\
\TAB \TAB 	$\}$; \\
\TAB \texttt{if} ($J+L \neq (J+L : x_n^\infty)$ \texttt{OR} $L = (1)$) \texttt{then} $\op{TerminationCondition}=$ \texttt{true}; \\
\TAB $J = J+L$;\\
\TAB \} \\
\texttt{return} $C+L$; \\
\}

\proof[Proof of correctness]
By Proposition~\ref{Prop:satulemmahard} at each step of the \texttt{while} loop, the ideal $L$ is such that $\V(L)$ contains the set of $a$ such that $J_a \neq (J_a:x_n)$. Moreover, by Noetherianity and by Proposition~\ref{Prop:satulemmahard} the \texttt{while} loop terminates after finitely many steps. When the algorithm terminates we have $J+L \neq (J+L : x_n^\infty)$, in which case $\V(L) = \{a \mid J_a \neq (J_a:x_n)\}$, or $L = (1)$, in which case $\V(L) = \{a \mid J_a \neq (J_a:x_n)\} = \emptyset$.
\endproof
\end{algoritmo}

Algorithm~\ref{satualgo} can be adapted to describe the set of parameters $a$ such that $J_a$ is not saturated with respect to $x_0\cdots x_n$.

\begin {algoritmo}\label{satualgon}
\leavevmode \\
{\bf {Input:}}  
			\begin{itemize}
				\item[$J$,] ideal in $\IC[c_1,\ldots, c_k, x_0, \ldots, x_n]$ homogeneous with respect to the $x_i$'s and saturated with respect to $x_0\cdots x_n$.
			\end{itemize}
{\bf {Output:}}  
			\begin{itemize}
				\item[$L$,] ideal in $\IC[c_1, \ldots, c_k]$ satisfying $\V(L) = \{a \mid J_a \neq (J_a:(x_0\cdots x_n))\}$.
			\end{itemize}
\{ \\
\texttt{for} $i$ \texttt{from} $0$ \texttt{to} $n$ \texttt{do} \\
\TAB \{ \\
\TAB let $J_i$ be the ideal obtained from $J$ by swapping $x_n$ and $x_i$; \\
\TAB let $L_i$ be the output of Algorithm~\ref{satualgo} on input $J_i$;\\
\TAB \}\\
$L=\sum L_i$;\\
\texttt{return} $L$;\\
\}
\proof[Proof of correctness]
As $(J_a: x_1 \cdots x_n) = (\ldots ((J_a:x_1):x_2)\ldots :x_n)$, we have that $J_a \neq (J_a:(x_1\cdots x_n))$ if and only if $J_a \neq (J_a:x_i)$ for some $i$. The correctness now follows from the correctness of Algorithm~\ref{satualgo}.
\end{algoritmo}

We are now ready to describe the main algorithm of this section. For the sake of clarity, given an ideal $L\subset \IC[c_1, \ldots, c_k]$ we will denote by $\dim_{\IA^k}(L)$ its dimension as an ideal of $\IC[c_1, \ldots, c_k]$, as opposed to its dimension as an ideal of $\IC[c_1, \ldots, c_k, x_1, \ldots, x_n]$.

\begin {algoritmo}\label{dimalgo}
\leavevmode \\
{\bf {Input:}}  
			\begin{itemize}
				\item[$I$,] ideal in $\IC[c_1,\ldots, c_k, x_1, \ldots, x_n]$ saturated with respect to $x_1 \cdots x_n$.
			\end{itemize}
{\bf {Output:}}  
			\begin{itemize}
				\item[$L$,] ideal in $\IC[c_1, \ldots, c_k]$ satisfying $\V(L) = \overline{\{a \mid 0 \leq \dim(\Va(I_a)) \leq 1\}}$.
			\end{itemize}
\{ \\
$C= I \cap \IC[c_1, \ldots, c_k]$;\\
$\{C_1, \ldots, C_m\} =$ \texttt{PrimaryDecomposition}($C$);\\
let $\bar{I}$ be the homogenization of $I$ in $(\IC[c_1, \ldots, c_k])[x_0, \ldots, x_n]$;\\
\texttt{for} $i$ \texttt{from} $1$ \texttt{to} $m$ \texttt{do} \\
\TAB \{ \\
\TAB $I_i = \bar{I} + C_i$; \\
\TAB $C' = C_i$; \\
\TAB $I' = I_i$; \\
\TAB $L_i = (0) \subseteq \IC[c_1, \ldots, c_k]$;\\
\TAB \texttt{while} ($\dim_{\IA^k}(C') + 1 \leq  \dim (I')$) \texttt{do} \\
\TAB \TAB \{ \\
\TAB \TAB $\{Q_{i1}, \ldots, Q_{il} \} =$ \texttt{PrimaryDecomposition}($I'$); \\
\TAB \TAB \texttt{for} $j$ \texttt{from} $1$ \texttt{to} $l$ \texttt{do} \\
\TAB \TAB \TAB \{ \\
\TAB \TAB \TAB $D_j = Q_{ij} \cap \IC[c_1, \ldots, c_k]$; \\
\TAB \TAB \TAB \texttt{if} ($\op{Radical}(D_j) = \op{Radical}(C_i)$ \texttt{AND} $\dim(Q_{ij}) > \dim_{\IA^k} (C_i) + 1 $) \texttt{then}  \\
\TAB \TAB \TAB \TAB \{compute the output $L_{ij}$ of Algorithm~\ref{satualgon} on input $Q_{ij} + \op{Radical}(C_i)$\}\\
\TAB \TAB \TAB \texttt{else} $L_{ij}=(0) \subseteq \IC[c_1, \ldots c_k]$; \\
\TAB \TAB \TAB \} \\
\TAB \TAB $L_i = \sum_j L_{ij}$;  \\
\TAB \TAB $I' = (I'+L_i:(x_0\cdots x_n)^\infty)$; \\
\TAB \TAB $C' = I' \cap \IC[c_1, \ldots, c_k]$;\\
\TAB \TAB $L_i = C'+L_i$;\\
\TAB \TAB \} \\
\TAB \} \\
$L= C +  \bigcap_i L_i$; \\
\texttt{return} $L$; \\
\}

\proof[Proof of correctness]
Let $\{C_1, \ldots, C_m\}$ be the primary ideals of a primary decomposition of $C$. The output of the Algorithm is the ideal $L= C +  \bigcap_i L_i$, so it will suffices to show that $\V(C_i + L_i)$ is the closure of the set  $\{a\in \V(C_i) \mid 0 \leq \dim(\Va(I_a)) \leq 1\}$ for $i = 1, \ldots, m$. 

  Let $C_i$ be a primary component of $C$ such that $\dim_{\IA^k} (C_i) +1 \geq \dim (I+C_i)$. The Algorithm computes $L_i = (0)$. For $a \in \V(C_i)$ general we have $0\leq \dim (I_a) \leq 1$. Therefore we have, as claimed, $\V(C_i + L_i) = \V(C) = \overline{\{a\in \V(C_i) \mid 0 \leq \dim(\Va(I_a)) \leq 1\}}$.  
	
	Now let $C_i$ be a primary component of $C$ such that $\dim_{\IA^k} (C_i) +1 < \dim (I+C_i)$. 
	We first show that the \texttt{while} loop terminates after finitely many steps. To do so, it suffices to show that the ideal $L_i$ as computed at the $j$\textsuperscript{th} iteration of the loop strictly contains the deal $L_i$ as computed at the $(j-1)$\textsuperscript{th} iteration. At the end of each iteration of the \texttt{while} loop the ideal $I'$ is redefined to contain $L_i$. Therefore it suffices to show that, at each iteration, the ideal $L_i$ is not contained in $I'$. Since $\dim_{\IA^k} (C_i) +1 < \dim (I')$ there exists some primary component $Q_{ij}$ of $I'$ such that $\V(Q_{ij})$ projects dominantly to $\V(C_i)$ and $\dim (C_i)+1 < \dim (Q_{ij})$. This matches the condition of the inner \texttt{if} operator, therefore the Algorithm computes the output $L_{ij}$ of Algorithm~\ref{satualgon}. By the correctness of Algorithm~\ref{satualgon} and Lemma~\ref{satulemma1} we have that $L_{ij}$, and hence $L_i$, is not contained in $I'$. This shows that the \texttt{while} loop terminates.
	We now show that the set of parameters $a\in \V(C_i)$ such that $0\leq \dim (I_a) \leq 1$ is contained in $V(C_i+L_i)$ at each iteration of the \texttt{while} loop. As before, let $Q_{ij}$ be a primary component of $I'$ such that $\V(Q_{ij})$ projects dominantly to $\V(C_i)$ and $\dim _{\IA^k}(C_i)+1 < \dim (Q_{ij})$.	By~\cite[Theorem~2 of Chapter~1 Section~8]{Mumred} the fibers of the projective closure of $\V(Q_{ij})$ all have dimension at least $2$. In particular the fiber $\Va((Q_{ij})_a)$ has dimension less than $2$ if and only if $\Va(Q_{ij})$ is empty. As a result the set of parameters $a$ such that $\Va((Q_{ij})_a)$ is empty is contained in the set of parameters $a$ such that $(Q_{ij})_a$ is not saturated with respect to $x_0\cdots x_n$. Since $\Va(I_a)$ contains $\Va((Q_{ij})_a)$, we have that the set of parameters $a\in \V(C_i)$ such that $0\leq \dim (I_a) \leq 1$ is contained in $\V(C_i+L_i)$. 
	Finally, since the \texttt{while} loop terminates, we have at the last iteration $\dim_{\IA^k} (C_i) +1 \geq \dim (I+C_i)$. Therefore the proof that $\V(C_i + L_i) = \overline{\{a\in \V(C_i) \mid 0 \leq \dim(\Va(I_a)) \leq 1\}}$ can be reduced to the previous case.
	\endproof
\end{algoritmo}

	\subsection{Hilbert functions}\label{Ssec:HiHil}

Theorem~\ref{prop:keymult} is a useful tool to compute the multiplicity of $\pos(e_1)$ on the fibers of $\V(I)$. We now discuss how to use this result to compute the locus in $\V(C)$ where the multiplicity is at least $m$.

When the projection $\pi:\V(J)\rightarrow \V(C)$ is quasi finite, the degree of $\pi$ equals the degree of the general fiber $\V(J_a)$ and hence, by Proposition~\ref{prop:torusfibers} and Theorem~\ref{prop:keymult}, the multiplicity of $\op{pos}(e_1)$ in $\Trop(I_a)$.  This gives, however, no control on the special fibers, as the following example shows. 

\begin{figure}
\begin{center}
\begin{minipage}{0.3\textwidth}
	\begin{center}
		\begin{tikzpicture}[style=thick]
			\draw [fill=none] (-1,1.5cm) circle(0pt) node[above]{$\Sigma_1$};
			\draw [->](0,0cm) -- (0,1cm);
			\draw [->](0,0cm) -- (-1,-1cm);
			\draw [->](0,0cm) -- (1,0cm);
			\draw [fill=none](0,1cm) circle (0pt) node[above] {$(0,1)$};
			\draw [fill=none](-1,-1cm) circle (0pt) node[below left] {$(-1,-1)$};
			\draw [fill=none](1,0cm) circle (0pt) node[right] {$(1,0)$};
			\draw [fill=white](0.5,0cm) circle (6pt) node[] {$1$};
			\draw [fill=white](0,0.5cm) circle (6pt) node[] {$1$};
			\draw [fill=white](-0.5,-0.5cm) circle (6pt) node[] {$1$};
		\end{tikzpicture}
	\end{center}
\end{minipage}
\begin{minipage}{0.3\textwidth}
	\begin{center}
		\begin{tikzpicture}[style=thick]
			\draw [fill=none] (-1,1.5cm) circle(0pt) node[above]{$\Sigma_2$};
			\draw [->](0,0cm) -- (0,1cm);
			\draw [->](0,0cm) -- (-1,-1cm);
			\draw [->](0,0cm) -- (1,0cm);
			\draw [fill=none](0,1cm) circle (0pt) node[above] {$(0,1)$};
			\draw [fill=none](-1,-1cm) circle (0pt) node[below left] {$(-1,-1)$};
			\draw [fill=none](1,0cm) circle (0pt) node[right] {$(1,0)$};
			\draw [fill=white](0.5,0cm) circle (6pt) node[] {$2$};
			\draw [fill=white](0,0.5cm) circle (6pt) node[] {$2$};
			\draw [fill=white](-0.5,-0.5cm) circle (6pt) node[] {$2$};
		\end{tikzpicture}
	\end{center}
\end{minipage}
\caption{The tropical varieties of Example~\ref{eg:sindown}}
\label{fig:sindown}
\end{center}
\end{figure}
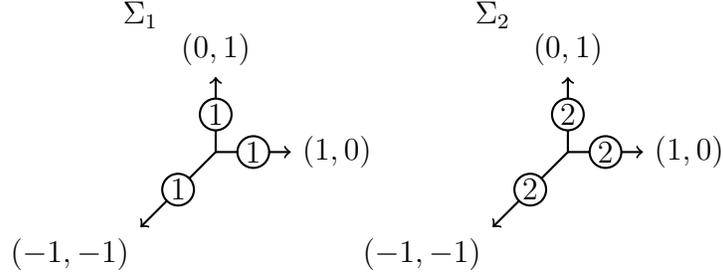

\begin{example}\label{eg:sindown}
Let $I= (c_1^3+c_1^2-c_2^2,2xc_2+2yc_2-c_1,2xc_1^2+2yc_1^2+2xc_1+2yc_1-c_2,4x^2c_1+8xyc_1+4y^2c_1+4x^2+8xy+4y^2-1) \subseteq \IC[c_1,c_2][x,y]$. The tropical variety $\Trop(I_a)$ is empty for $a\notin \V(c_1^3+c_1^2-c_2^2)$, it is the tropical curve $\Sigma_1$ depicted in Figure~\ref{fig:sindown} for $(0,0) \neq a\in \V(c_1^3+c_1^2-c_2^2)$, and it is  the tropical curve $\Sigma_2$ depicted in Figure~\ref{fig:sindown} for $a=(0,0)$. 
\end{example}

Proposition~\ref{prop:keymult} relates the multiplicity of $e_1$ in $\Trop(I_a)$ with the cardinality of $\overline{\V(I_a)} \cap \BO_1$. The following proposition relates this cardinality with the Hilbert function of $J_a$.



\begin{proposition}\label{prop:hilbfun}
Suppose that for $a\in \V(C)$ general the dimension of $\V(J_a)$ is $0$. Denote by $\bar{J}$ the homogenization of $J$ in $(\IC[c_1, \ldots, c_k])[x_0, \ldots, x_n]$, and by $\bar{J}_a$ the homogenization of $J_a$ in $\IC[x_0, \ldots, x_n]$.
Let $r$ be the Castelnuovo-Mumford regularity $\op{reg}(\bar{J})$ and denote by $h_{\bar{J}_a}$ the Hilbert polynomial of $\bar{J}_a$.
Then for every  $a\in \V(C)$ we have $\mult(\op{pos}(e_1),I_a) \leq h_{\bar{J}_a}(r)$, and the equality holds on an open dense subset of  $\V(C)$.

\proof
We recall that $\RReal_1$ is defined as 
\begin {equation*}
	\RReal_1 \eeq \{ (a,x) \in \IA^{k+n} \mid x\in \overline{\Va(I_a)} \cap \BO_1\}.
\end{equation*}
Fix $a\in \V(C)$ such that $\dim(\Va(I_a)) \leq 1$. The multiplicity $\mult(\op{pos}(e_1),I_a) \leq h_{\bar{J}_a}(r)$ equals the cardinality of $(\RReal_1)_a$. By Proposition~\ref{prop:torusfibers}, $(\RReal_1)_a \subseteq \V(J_a)$ for every $a\in \V(C)$ and, by Theorem~\ref{localgteo2}, $(\RReal_1)_a = \V(J_a)$ for $a\in \V(C)$ general. Moreover, since $\dim(J_a)$ is $0$, its cardinality equals the cardinality of its closure $\V(\bar{J}_a)$ in the projective space. To conclude the proof it suffices to show that the Hilbert function of $\bar{J}_a$ computed at $r$ coincides with the Hilbert polynomial of $\bar{J}_a$ computed at $r$. This follows from \cite[Corollary~20.19]{Eis} once we identify the ideal $\bar{J}_a\subseteq \IC[x_0, \ldots, x_n]$ with the ideal $(\bar{J}, c_1-a_1, \ldots, c_k-a_k) \subseteq \IC[c_1, \ldots, c_k, x_0,\ldots, x_n]$.
\endproof
\end{proposition}

We now describe how to compute the set of parameters $a\in \V(C)$ such that $h_{\bar{J}_a}(r)\geq m$. 
We denote by $\bar{J}(r)$ the vector space $\bar{J} \cap  (\IC[c_1, \ldots, c_k])[ x_0, \ldots, x_n]_r$.
Let $\BL=\{f_1, \ldots, f_l\}$ be a homogeneous basis of $\bar{J}$. A set of generators of $\bar{J}(r)$ as a vector space is given by all the degree $r$ polynomials that can be obtained by multiplying an element of $\BL$ by a monomial. Let $M$ be the coefficient matrix of this basis. It is a matrix with entries in $\IC[c_1, \ldots, c_k]$. Let $D$ the ideal generated by $\left(\binom{n+r}{r}-m+1\right)$-minors of $M$. For $a\in \V(C)$ the degree $r$ component of $\bar{J}_a$ is a  vector space $(\bar{J}_a)(r) \subseteq   \IC[ x_0, \ldots, x_n]_r$.  Denote by $M_a$ the matrix obtained by evaluating the entries of $M$ in $a$. By construction $M_a$ is the coefficients matrix of a set of generators of $\bar{J}_a(r)$. In particular $\bar{J}_a(r)$ has dimension  at most $\binom{n+r}{r}-m$ if and only if $a\in \V(D)$ and, equivalently,  we have $h_{\bar{J}_a}(r)\geq m$ if and only if $a\in \V(D)$. 
This construction is the content of the following algorithm. 

\begin{algoritmo}\label{Algo:macaulayminors}
\leavevmode \\
{\bf {Input:}}  
			\begin{itemize}
				\item[$I$,] ideal in $\IC[c_1,\ldots, c_k, x_1, \ldots, x_n]$ with the property that, for $a\in \V(I)\cap \IA^k$ general, $\dim (\Va(I_a)) \leq 1$,
				\item[$m$,] positive integer.
			\end{itemize}
{\bf {Output:}}  
			\begin{itemize}
				\item[$D$,] ideal in $\IC[c_1, \ldots, c_k]$ such that $h_{\bar{J}_a}(r)\geq m$ if and only if $a\in \V(D)$.
			\end{itemize}
\{ \\
let $E$ be the homogenization of $J$ in $(\IC[c_1,\ldots, c_k]) [x_0,x_1, \ldots, x_n]$ and $L$ a list of generators of $E$; \\
$B= \{\}$; \\
\texttt{for} $f$ \texttt{in} $L$ \texttt {do}  \\
\TAB \{ add every monomial $x^u$ of degree $m-\op{deg}(f)$ to $B$ \} ; \\
let $M$ be the matrix of coefficients of $B$; \\
compute the ideal $D$ generated by the $(\binom{n+r}{r}-m+1)$-minors of $M$; \\
\texttt{return} $D$; \\
\} \\
\end{algoritmo}

\begin{proposition}\label{prop:minorsideal}
Suppose that, for $a\in \V(C)$ general, we have $\dim (\Va(I_a)) \leq 1$. Let $a\in \V(C)$  be a parameter such that $ \mult(\op{pos}(e_1),\Trop(I_a))\geq m$, then $a\in\V(D)$. In other words $\Real_{1,m} \subseteq \V(C, D)$.
\proof
As $ \mult(\op{pos}(e_1),\Trop(I_a))\geq m$ then, by Proposition~\ref{prop:hilbfun}, $h_{J_a}(r) \geq m$. As a result the matrix $M_a$ has rank at most $\binom{n+r}{r}-m$ and therefore $a\in\V(D)$.
\endproof
\end{proposition}


When $D \subseteq C$ the ideal of  $\Real_{1,m}$ is described by the following proposition.

\begin{proposition}\label{prop:termcond}
Suppose that, for $a\in \V(C)$ general, we have $\dim (\Va(I_a)) \leq 1$ and suppose that $D\subseteq C$. Then we have $ \I(\Real_{1,m}) =C$.
\proof
By Proposition~\ref{prop:minorsideal} we have $ \I(\Real_{1,m}) \supset C+D = C$.
Moreover by the definition of $D$ we have $h_{J_a}(r) \geq m$ for every $a \in \V(C,D)= \V(C)$ and by Proposition~\ref{prop:hilbfun} we have that  $h_{J_a}(r) = \mult(\op{pos}(e_1),\Trop(I_a))$ for $a$ general in $\V(C)$. This shows that  $\Real_{1,m}$ is dense in $\V(C)$.
\endproof
\end{proposition}

	\subsection{The main algorithm}\label{Ssec:HiGlu}

We are now ready to describe an algorithm to compute the Zariski closure of $\Real_{1,m}$. The algorithm proceeds as follows. It first computes the ideals $C$ and $J$ using Algorithm~\ref{LOCALG} on input $I$. It then computes the ideals $J$, $C$ and $D$. If $D\not\subseteq C$ it starts again running Algorithm~\ref{LOCALG} on the restricted family $I+C+D$.

\begin{algoritmo}\label{weightalgo}
\leavevmode \\
{\bf {Input:}}  
			\begin{itemize}
				\item[$I$,] ideal in $\IC[c_1,\ldots, c_k, x_1, \ldots, x_n]$,
				\item[$m$,] positive integer.
			\end{itemize}
{\bf {Output:}}  
			\begin{itemize}
				\item[$D$,] ideal in $\IC[c_1, \ldots, c_k]$ satisfying $\V(D)=\overline{\Real}_{1,m}\subseteq \IA^k$.
			\end{itemize}
\{ \\
$C = D = (0) \subseteq \IC[c_1, \ldots, c_k]$; \\
$\op{doFirstIteration} =$ \texttt{true}; \\
\texttt{while} ($D \not\subseteq C$ \texttt{ OR } $\op{doFirstIteration} =$ \texttt{true})  \texttt{do} \\
\TAB $\{$ \\
\TAB 	 $\op{doFirstIteration} =$ \texttt{false}; \\
\TAB 	 $\op{computeLocCond} =$ \texttt{true}; \\
\TAB 	 \texttt{while}($\op{computeLocCond}$) \texttt{do} \\
\TAB \TAB $\{$ \\
\TAB \TAB compute the output $(C,J)$ of Algorithm~\ref{LOCALG} with input $I+C+D$; \\
\TAB \TAB compute the output $C'$ of Algorithm~\ref{dimalgo} with input $(I+C:(x_1\cdots x_n)^\infty)$;\\ 
\TAB \TAB \texttt{if} $C =  C'$ \texttt{then} $\op{computeLocCond} =$ \texttt{false} ; \\
\TAB \TAB $C=C'$; \\
\TAB \TAB \}; \\
\TAB 	 compute the output $D$ of Algorithm~\ref{Algo:macaulayminors} on input $(I+C, m)$; \\
\TAB 	 $\}$; \\
\texttt{return} $C$; \\
\} \\

\proof[Proof of correctness]
When the algorithm reaches its termination we are in assumptions of Proposition~\ref{prop:termcond} and therefore $\overline{\Real}_{1,m} =\V(C)$.
\endproof
\end{algoritmo}

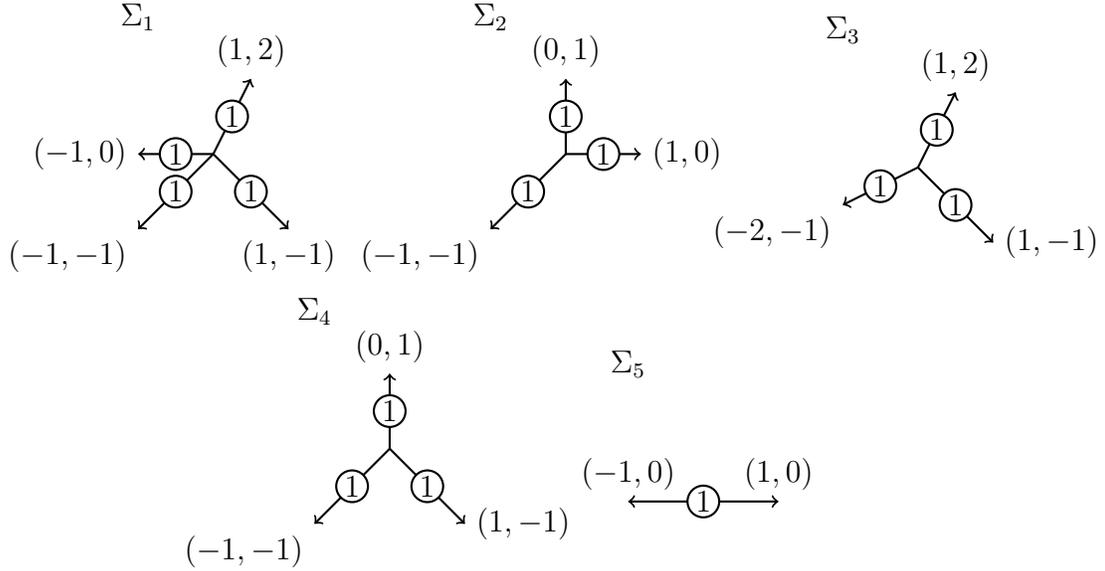
\begin{figure}
\begin{center}
\begin{minipage}{0.3\textwidth}
	\begin{center}
		\begin{tikzpicture}[style=thick]
			\draw [fill=none] (-1,1.5cm) circle(0pt) node[above]{$\Sigma_1$};
			\draw [->](0,0cm) -- (0.5,1cm);
			\draw [->](0,0cm) -- (-1,-1cm);
			\draw [->](0,0cm) -- (-1,0cm);
			\draw [->](0,0cm) -- (1,-1cm);
			\draw [fill=none](0.5,1cm) circle (0pt) node[above] {$(1,2)$};
			\draw [fill=none](-1,-1cm) circle (0pt) node[below left] {$(-1,-1)$};
			\draw [fill=none](-1,0cm) circle (0pt) node[left] {$(-1,0)$};
			\draw [fill=none](1,-1cm) circle (0pt) node[below] {$(1,-1)$};
			\draw [fill=white](-0.5,0cm) circle (6pt) node[] {$1$};
			\draw [fill=white](0.25,0.5cm) circle (6pt) node[] {$1$};
			\draw [fill=white](-0.5,-0.5cm) circle (6pt) node[] {$1$};
			\draw [fill=white](0.5,-0.5cm) circle (6pt) node[] {$1$};
		\end{tikzpicture}
	\end{center}
\end{minipage}
\begin{minipage}{0.3\textwidth}
	\begin{center}
		\begin{tikzpicture}[style=thick]
			\draw [fill=none] (-1,1.5cm) circle(0pt) node[above]{$\Sigma_2$};
			\draw [->](0,0cm) -- (0,1cm);
			\draw [->](0,0cm) -- (-1,-1cm);
			\draw [->](0,0cm) -- (1,0cm);
			\draw [fill=none](0,1cm) circle (0pt) node[above] {$(0,1)$};
			\draw [fill=none](-1,-1cm) circle (0pt) node[below left] {$(-1,-1)$};
			\draw [fill=none](1,0cm) circle (0pt) node[right] {$(1,0)$};
			\draw [fill=white](0.5,0cm) circle (6pt) node[] {$1$};
			\draw [fill=white](0,0.5cm) circle (6pt) node[] {$1$};
			\draw [fill=white](-0.5,-0.5cm) circle (6pt) node[] {$1$};
		\end{tikzpicture}
	\end{center}
\end{minipage}
\begin{minipage}{0.3\textwidth}
	\begin{center}
		\begin{tikzpicture}[style=thick]
			\draw [fill=none] (-1,1.5cm) circle(0pt) node[above]{$\Sigma_3$};
			\draw [->](0,0cm) -- (0.5,1cm);
			\draw [->](0,0cm) -- (-1,-0.5cm);
			\draw [->](0,0cm) -- (1,-1cm);
			\draw [fill=none](0.5,1cm) circle (0pt) node[above] {$(1,2)$};
			\draw [fill=none](-1,-0.5cm) circle (0pt) node[below left] {$(-2,-1)$};
			\draw [fill=none](1,-1cm) circle (0pt) node[right] {$(1,-1)$};
			\draw [fill=white](0.5,-0.5cm) circle (6pt) node[] {$1$};
			\draw [fill=white](0.25,0.5cm) circle (6pt) node[] {$1$};
			\draw [fill=white](-0.5,-0.25cm) circle (6pt) node[] {$1$};
		\end{tikzpicture}
	\end{center}
\end{minipage}
\begin{minipage}{0.3\textwidth}
	\begin{center}
		\begin{tikzpicture}[style=thick]
			\draw [fill=none] (-1,1.5cm) circle(0pt) node[above]{$\Sigma_4$};
			\draw [->](0,0cm) -- (0,1cm);
			\draw [->](0,0cm) -- (-1,-1cm);
			\draw [->](0,0cm) -- (1,-1cm);
			\draw [fill=none](0,1cm) circle (0pt) node[above] {$(0,1)$};
			\draw [fill=none](-1,-1cm) circle (0pt) node[below left] {$(-1,-1)$};
			\draw [fill=none](1,-1cm) circle (0pt) node[right] {$(1,-1)$};
			\draw [fill=white](0.5,-0.5cm) circle (6pt) node[] {$1$};
			\draw [fill=white](0,0.5cm) circle (6pt) node[] {$1$};
			\draw [fill=white](-0.5,-0.5cm) circle (6pt) node[] {$1$};
		\end{tikzpicture}
	\end{center}
\end{minipage}
\begin{minipage}{0.3\textwidth}
	\begin{center}
		\begin{tikzpicture}[style=thick]
			\draw [fill=none] (-1,1.5cm) circle(0pt) node[above]{$\Sigma_5$};
			\draw [->](0,0cm) -- (1,0cm);
			\draw [->](0,0cm) -- (-1,0cm);
			\draw [fill=none](1,0cm) circle (0pt) node[above] {$(1,0)$};
			\draw [fill=none](-1,0cm) circle (0pt) node[above] {$(-1,0)$};
			\draw [fill=white](0,0cm) circle (6pt) node[] {$1$};
		\end{tikzpicture}
	\end{center}
\end{minipage}
\caption{The tropical varieties of Example~\ref{eg5}}
\label{fig:eg5}
\end{center}
\end{figure}

\begin{example}\label{eg5}
Let $I$ be the ideal $I = (c_1y+ xy+xy^2+c_1c_2x^2 + c_2x^2y + x^2y)\subseteq \IC[c_1, c_2, x, y]$. The tropical variety $\Trop(I_a)$ is depicted in Figure~\ref{fig:eg5}: 
it is the tropical curve $\Sigma_1$ for $a_1 \neq 0, 1_2\neq 0,-1$, 
it is the tropical curve $\Sigma_2$ for $a_1= 0, a_2 \neq -1$, 
it is the tropical curve $\Sigma_3$ for $a_1 \neq 0, a_2 = -1$, 
it is the tropical curve $\Sigma_4$ for $a_1 \neq 0, a_2 = 0$ 
and it is the tropical curve $\Sigma_5$ for $a_1= 0, a_2 = -1$. 
Let $m=2$, we have $\Real_{1,2} = \emptyset$. Algorithm~\ref{weightalgo} first computes the ideals $C= \I(\Real_1) = (c_1)$ and $J= \I(\RReal_1) = (c_1, x, y+1)$. Finally it computes that the Hilbert function of $I_0$ is less than $2$.
\end{example}

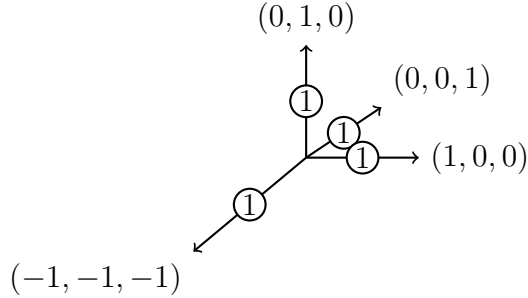
\begin{figure}
\begin{center}
\begin{minipage}{0.5\textwidth}
	\begin{center}
		\begin{tikzpicture}[style=thick]
			\draw [->](0,0cm) -- (0,1.5cm);
			\draw [->](0,0cm) -- (-1.5,-1.25cm);
			\draw [->](0,0cm) -- (1.5,0cm);
			\draw [->](0,0cm) -- (1,0.67cm);
			\draw [fill=none](0,1.5cm) circle (0pt) node[above] {$(0,1,0)$};
			\draw [fill=none](-1.5,-1.25cm) circle (0pt) node[below left] {$(-1,-1,-1)$};
			\draw [fill=none](1.5,0cm) circle (0pt) node[right] {$(1,0,0)$};
			\draw [fill=none](1,0.67cm) circle (0pt) node[above right] {$(0,0,1)$};
			\draw [fill=white](0.5,0.33cm) circle (6pt) node[] {$1$};
			\draw [fill=white](-0.75,-0.625cm) circle (6pt) node[] {$1$};
			\draw [fill=white](0.75,0cm) circle (6pt) node[] {$1$};
			\draw [fill=white](0,0.75cm) circle (6pt) node[] {$1$};
		\end{tikzpicture}
	\end{center}
\end{minipage}
\caption{The tropical variety $\Trop(\V(I_a))$ of Example~\ref{eg9}.}
\label{fig:eg9}
\end{center}
\end{figure}

\begin{example}\label{eg9}
Let $I$ be the ideal $I= I_0 \cdot I_1 \subseteq \IC[c,x,y,z]$ where $I_0=(x+y+z+1, 2x-3y+5z-2, c)$ and $I_1=(x+3y-2z+6, c(-x+4y+2z-1))$. The tropical variety $\Trop(I_a)$ does not depend on $a$ and it is depicted in Figure~\ref{fig:eg9}. For $m=2$ we have $\Real_{1,2} = \emptyset$. The algorithm computes the ideal $\I(\Real_{1,2}) = \IC[c]$ in the following way. It first computes $C= \I(\Real_1(I)) = (0) \subseteq \IC[c]$ and $J= \I(\RReal_1(I))= I + (x)$. Then it computes that the Hilbert function $h_{J_a} (2)$ of $J_a$ computed at the regularity $2$ is at least $2$ in $\V(c)$. Then the algorithm computes $C= \I(\Real_1(I+(c))) = (c) \subseteq \IC[c]$ and $J= \I(\RReal_1(I+(c)))= I + (x,c)$. The ideal $J$ has two components, corresponding to the two components $I_0$ and $I_1$ of $I$, and the component corresponding to $I_1$ is discarded since the fiber of $I_1+(c)$ has dimension $2$. Finally, the computation of the Hilbert function shows that $\I(\Real_{1,2}) = \IC[c]$.
\end{example}

\section{Change of Coordinates}\label{sec4}

In the previous sections we studied the realizability space $\Real_{1,m}$ for the ray $\pos(e_1)$ with multiplicity $m$. In this section we show how to extend these results to the case of an arbitrary rational ray $\rho$. We denote by $\Real_{\rho, m}$ the set of parameters $a\in \IA^k$ such that the tropicalization of $X_a$ is a curve that contains the ray $\rho$ with multiplicity at least $m$.

To a $m\times n$ integer matrix $A=(a_{i,j})$ one can associate the monomial map $\varphi_A$, that is the regular morphism 
\begin{equation*}
\begin{array}{cccc}
\varphi_A: & \T^n			 					& \rightarrow & \T^m \\
				   & (x_1, \ldots, x_n) & \mapsto			& (\prod_{i=1}^n x_i^{a_{1,i}}, \ldots, \prod_{i=1}^n x_i^{a_{m,i}}).
\end{array}
\end{equation*}
The tropicalization $\Trop(\varphi_A)$ of $\varphi$ is by definition the linear map $\Trop(\varphi_A): \IR^n\rightarrow \IR^m$ associated to the transpose matrix $A^\top$.
The reason for this definition is the following (see \cite[Corollary 3.2.15]{Mac}): for any subvariety $X$ of $\T^n$  the tropicalization of the Zariski closure of its image under $\varphi_A$ equals the image of its tropicalization under $\Trop(\varphi_A)$, in symbols
\begin{equation}\label{eq:tropfunct}
		\Trop(\overline{\varphi_A(X)})=\Trop(\varphi_A)(\Trop(X)).
\end{equation}

Let $\rho\subseteq \IR^n$ be a rational ray and let $v$ be the first integer point of $\rho$. By the Lemma in~\cite{Schenkman} there exists a $\op{GL}_n(\IZ)$ matrix $B$ such that $B\cdot e_1 = v$. We denote by $A$ the transpose of $B$. Let $\varphi_A:T^n\rightarrow T^n$ be the monomial map associated to $A$, and denote by $\varphi_A^*:\IC[x_1^\pm, \ldots, x_n^\pm]\rightarrow \IC[x_1^\pm, \ldots, x_n^\pm]$ the associated ring homomorphism. As $A\in \op{GL}_n(\IZ)$ is invertible, $\varphi_A^*$ is invertible with inverse $\varphi_{A^{-1}}^*$. Given a variety $X\subseteq \T^n$, by Equation~\eqref{eq:tropfunct} we have that $\mult(\op{pos}(v),\Trop(X)) = \mult(\op{pos}(e_1),\Trop(\overline{\varphi_A(X)}))$. The variety $\overline{\varphi_A(X)}$ is easily described: if $X$ is the vanishing locus of the Laurent polynomials $f_1, \ldots, f_l \in \IC[x_1^\pm, \ldots, x_n^\pm]$, then $\overline{\varphi_A(X)}$ is the vanishing locus of the Laurent polynomials $(\varphi_A^*)^{-1}(f_1), \ldots, (\varphi_A^*)^{-1}(f_l)$.

Let now $I \subseteq \IC[c_1, \ldots, c_k][x_1, \ldots, x_n]$. We denote again by $\varphi_A^*$ the ring homomorphism $\varphi_A^*:\IC[c_1, \ldots, c_k][x_1^\pm, \ldots, x_n^\pm]\rightarrow \IC[c_1, \ldots, c_k][x_1^\pm, \ldots, x_n^\pm]$. We have the following equality:
\begin{equation*}
	\Real_{\rho,m} (I) = \Real_{1,m} ((\varphi_A^*)^{-1} \ I).
\end{equation*}
This allows to compute $\Real_{\rho, m}(I)$ via Algorithm~\ref{weightalgo}. 

We conclude this section by describing an algorithm to explicitly compute a matrix $B \in \op{GL}_n(\IZ)$ satisfying the property $Be_1=v$.
The Smith normal form of an integer matrix $M\in \IZ^{n\times m}$ is a matrix $D\in \IZ^{n\times m}$ that is zero outside the main diagonal and whose elements on the diagonal $(d_1, \ldots, d_l)$ satisfy $d_1|d_2\ldots |d_l$, where $l=\min(n,m)$. The Smith normal form $D$ is related to $M$ via the expression $D=PMQ$ where $P$ is an invertible $n\times n$ integer matrix and $Q$ is an invertible $m\times m$ integer matrix.

\begin{algoritmo}\label{alg:fitsl}
\leavevmode \\
{\bf {Input:}}  
			\begin{itemize}
				\item[$v$,] a primitive integer vector $v\in \IZ^n$.
			\end{itemize}
{\bf {Output:}}  
			\begin{itemize}
				\item[$B$,] a $\op{GL}_n(\IZ)$ matrix such that $B e_1 = v$.
			\end{itemize}
\{\\
Compute the a Smith normal expression $D= Pv^\top Q$ for the row vector $v^\top $;\\
$d= D_{1,1}$;\\
$B = dP(Q^\top )^{-1}$;\\
\texttt{return} $B$;\\
\}\\

\proof[Proof of correctness]
The Smith normal form $D$ is a $1\times n$ integer matrix $D=(d,0,\ldots, 0)$. Moreover, $P$ is a $1\times 1$ invertible matrix, so we have $P=(1)$ or $P=(-1)$. We have $P  v^\top  Q = D= de_1^\top $ and hence $P  Q^\top   v = de_1$. It follows that $d P (Q^\top )^{-1} e_1 = d(P Q^\top )^{-1} e_1 = v$.
\endproof
\end{algoritmo}

\section{The Global Algorithm}\label{sec5}

In this section we glue the information coming from Algorithm~\ref{weightalgo} to describe the closure of the locus of parameters $a\in\IC^k$ such that $\Trop(I_a)$ is a curve and contains a given list of rays $\Sigma=\{\rho_1, \ldots, \rho_s\}$ with at least multiplicities $m_1, \ldots, m_s$.

\begin{algoritmo}\label{GLOBWEI}
\leavevmode \\
{\bf {Input:}}  
			\begin{itemize}
				\item[$I$,] ideal in $\IC[c_1,\ldots, c_k, x_1, \ldots, x_n]$,
				\item[$\Sigma$,] a collection of rays in $\IR^n$ and multiplicities $\Sigma=\{(\rho_1, m_1), \ldots, (\rho_s, m_s)\}$.
			\end{itemize}
{\bf {Output:}}  
			\begin{itemize}
				\item[$C$,] ideal in $\IC[c_1, \ldots, c_k]$ satisfying $\V(C)=\overline{\Real}_{\Sigma}\subseteq \IA^k$.
			\end{itemize}
\{ \\
$C=C_0=(0)\subseteq \IC[c_1,\ldots, c_k]$; \\
$\op{doFirstIteration} =$\texttt{true}; \\
\texttt{while} ($C \neq C_0$ or $\op{doFirstIteration}$) \texttt{do} \\
\TAB  $\{$ \\
\TAB  $\op{doFirstIteration} =$\texttt{false}; \\
\TAB	 $C_0=C$; \\
\TAB	\texttt{for} $(\rho_i, m_i)$ \texttt{in} $\Sigma$ \texttt{do} \\
\TAB	\TAB	 $\{$ \\
\TAB	\TAB	compute a matrix $A \in \op{GL}_n(\IZ)$ such that $A^\top \cdot v = e_1$ via Algorithm~\ref{alg:fitsl}; \\
\TAB	\TAB	compute $\I (\Real_{1,m_i}((\varphi_A^*)^{-1} (I) + C))$ via Algorithm~\ref{weightalgo};  \\
\TAB	\TAB	$C= \I( \Real_{1,m_i}((\varphi_A^*)^{-1} (I) + C))$; \\
\TAB	\TAB	 $\}$ \\
\TAB	 $\}$ \\
\texttt{return} $C$; \\
\}	 \\
\end{algoritmo}

We need to repeat the procedure in Algorithm~\ref{GLOBWEI} as Algorithm~\ref{weightalgo} assures that the multiplicity of $\op{pos}(e_1)$ in $\Trop(I_a)$ is at least $m$ only for a dense subset in its output. The following is an example of when a single iteration is not sufficient.

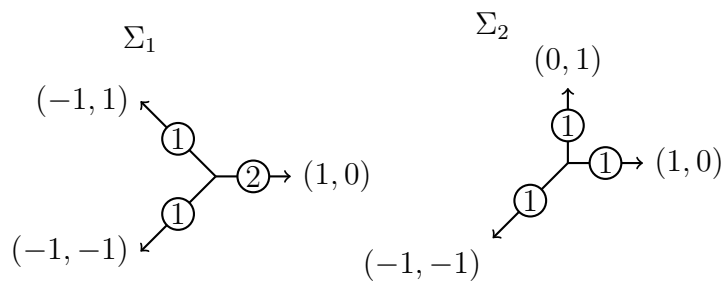
\begin{figure}
\begin{center}
\begin{minipage}{0.3\textwidth}
	\begin{center}
		\begin{tikzpicture}[style=thick]
			\draw [fill=none] (-1,1.5cm) circle(0pt) node[above]{$\Sigma_1$};
			\draw [->](0,0cm) -- (-1,1cm);
			\draw [->](0,0cm) -- (-1,-1cm);
			\draw [->](0,0cm) -- (1,0cm);
			\draw [fill=none](-1,1cm) circle (0pt) node[left] {$(-1,1)$};
			\draw [fill=none](-1,-1cm) circle (0pt) node[left] {$(-1,-1)$};
			\draw [fill=none](1,0cm) circle (0pt) node[right] {$(1,0)$};
			\draw [fill=white](0.5,0cm) circle (6pt) node[] {$2$};
			\draw [fill=white](-0.5,0.5cm) circle (6pt) node[] {$1$};
			\draw [fill=white](-0.5,-0.5cm) circle (6pt) node[] {$1$};
		\end{tikzpicture}
	\end{center}
\end{minipage}
\begin{minipage}{0.3\textwidth}
	\begin{center}
		\begin{tikzpicture}[style=thick]
			\draw [fill=none] (-1,1.5cm) circle(0pt) node[above]{$\Sigma_2$};
			\draw [->](0,0cm) -- (0,1cm);
			\draw [->](0,0cm) -- (-1,-1cm);
			\draw [->](0,0cm) -- (1,0cm);
			\draw [fill=none](0,1cm) circle (0pt) node[above] {$(0,1)$};
			\draw [fill=none](-1,-1cm) circle (0pt) node[below left] {$(-1,-1)$};
			\draw [fill=none](1,0cm) circle (0pt) node[right] {$(1,0)$};
			\draw [fill=white](0.5,0cm) circle (6pt) node[] {$1$};
			\draw [fill=white](0,0.5cm) circle (6pt) node[] {$1$};
			\draw [fill=white](-0.5,-0.5cm) circle (6pt) node[] {$1$};
		\end{tikzpicture}
	\end{center}
\end{minipage}
\caption{The tropical varieties of Example~\ref{eg:10}}
\label{fig:eg10}
\end{center}
\end{figure}

\begin{example}\label{eg:10}
Let $I$ be the ideal $I= (c+y+xy+y^2) \subseteq \IC[c, x,y]$. The tropical variety $\Trop(I_a)$ is depicted in Figure~\ref{fig:eg10}: it is $\Sigma_1$ for $a\neq 0$ and it is $\Sigma_2$ for $a=0$. Let $\Sigma = \{(\op{pos}(e_1), 2), (\op{pos}(e_2), 1)\}$. We have $\Real_\Sigma = \emptyset$. To compute the ideal $\I(\Real_\Sigma)$ Algorithm~\ref{GLOBWEI} needs two iterations of the \texttt{while} loop. It first computes $\Real_{\pos{e_1}, 2} (I)= \IA^1$ and $\Real_{\pos{e_2}, 1} (I)= \V(c)$. It then computes that $\Real_{\pos{e_1}, 2} (I+(c))= \emptyset$

\end{example}

\addcontentsline{toc}{section}{References}
\bibliographystyle{plain}
\bibliography{rrbib}

\def\cprime{$'$}
\begin{thebibliography}{10}

\bibitem{CLOS}
David Cox, John Little, and Donal O'Shea.
\newblock {\em Ideals, varieties, and algorithms}.
\newblock Undergraduate Texts in Mathematics. Springer-Verlag, New York, second
  edition, 1997.
\newblock An introduction to computational algebraic geometry and commutative
  algebra.

\bibitem{Eis}
David Eisenbud.
\newblock {\em Commutative algebra}, volume 150 of {\em Graduate Texts in
  Mathematics}.
\newblock Springer-Verlag, New York, 1995.
\newblock With a view toward algebraic geometry.

\bibitem{M2}
Daniel~R. Grayson and Michael~E. Stillman.
\newblock Macaulay2, a software system for research in algebraic geometry.
\newblock Available at \url{http://www.math.uiuc.edu/Macaulay2/}.

\bibitem{NL}
Phillip Griffiths and Joe Harris.
\newblock On the {N}oether-{L}efschetz theorem and some remarks on
  codimension-two cycles.
\newblock {\em Math. Ann.}, 271(1):31--51, 1985.

\bibitem{CompGBAlgo}
Deepak Kapur, Yao Sun, and Dingkang Wang.
\newblock An efficient method for computing comprehensive {G}r\"obner bases.
\newblock {\em J. Symbolic Comput.}, 52:124--142, 2013.

\bibitem{Mac}
Diane Maclagan and Bernd Sturmfels.
\newblock {\em Introduction to tropical geometry}, volume 161 of {\em Graduate
  Studies in Mathematics}.
\newblock American Mathematical Society, Providence, RI, 2015.

\bibitem{Mumred}
David Mumford.
\newblock {\em The red book of varieties and schemes}, volume 1358 of {\em
  Lecture Notes in Mathematics}.
\newblock Springer-Verlag, Berlin, expanded edition, 1999.
\newblock Includes the Michigan lectures (1974) on curves and their Jacobians,
  With contributions by Enrico Arbarello.

\bibitem{Schenkman}
Eugene Schenkman.
\newblock The basis theorem for finitely generated abelian groups.
\newblock {\em The American Mathematical Monthly}, 67(8):770--771, 1960.

\bibitem{Speyer}
David~E. Speyer.
\newblock {\em Tropical geometry}.
\newblock ProQuest LLC, Ann Arbor, MI, 2005.
\newblock Thesis (Ph.D.)--University of California, Berkeley.

\bibitem{BasesandPolytopes}
Bernd Sturmfels.
\newblock {\em Gr\"obner bases and convex polytopes}, volume~8 of {\em
  University Lecture Series}.
\newblock American Mathematical Society, Providence, RI, 1996.

\bibitem{Tev}
Jenia Tevelev.
\newblock Compactifications of subvarieties of tori.
\newblock {\em Amer. J. Math.}, 129(4):1087--1104, 2007.

\bibitem{ComprehensiveGB}
Volker Weispfenning.
\newblock Comprehensive {G}r\"obner bases.
\newblock {\em J. Symbolic Comput.}, 14(1):1--29, 1992.

\end{thebibliography}

\end {document}